\documentclass[12pt,leqno,a4paper]{article}
\usepackage{amsfonts}
\usepackage{amsmath, amssymb, amsthm, amsfonts}

\numberwithin{equation}{section}

\begin{document}
\title{On the error term in Weyl's law for the Heisenberg manifolds}
\author{ Wenguang ZHAI (Jinan)}
\date{}

\footnotetext[0]{2000 Mathematics Subject Classification: 11N37,
35P20, 58J50.} \footnotetext[0]{Key Words: Heisenberg manifold,
error term, exponential sum, power moment.} \footnotetext[0]{This
work is supported by  National Natural Science Foundation of
China(Grant No. 10301018) and  National Natural Science Foundation
of Shandong Province(Grant No. 2006A31).} \maketitle

{\bf Abstract.} For a fixed integer $l\geq 1$  , let $R(t)$ denote
the error term in the Weyl's law of  a $(2l+1)$-dimensional
Heisenberg manifold with the  metric $g_l.$ In this paper we shall
prove the asymptotic formula of the $k$-th power moment for any
integers $3\leq k\leq 9.$ We shall also prove that the function
$t^{-(l-1/4)}R(t)$ has a distribution function.

\bigskip

{\bf Contents}\\
$\S 1$. Introduction \\
1.1  The Weyl's law for ${\Bbb T}^2:$ the Gauss circle problem\\
1.2  The Weyl's law for $3$-dimensional Heisenberg manifold\\
$\S 2$.  Main Results\\
2.1 New results for power moments of $R(t)$\\
2.2 Idea of the proof\\
$\S 3$. Background of Heisenberg manifolds and the $\psi$-expression
of
$R(t)$\\
3.1 Heisenberg manifolds\\
3.2 The spectrum of Heisenberg manifolds\\
3.3 the $\psi$-expression of $R(t)$\\
3.4 A weighted lattice point problem\\
$\S 4$. Some preliminary Lemmas\\
$\S 5$. Proof of Theorem 1 \\
5.1 Large value estimate of $R_*(x)$\\
5.2 Proof of Theorem 1\\
$\S 6$. Proof of Theorem 2 \\
6.1 The upper bound of  $\int_T^{2T}G(x)dx$\\
6.2 An analogue of Voronoi's formula for $R_*(x)$\\
6.3 Evaluation of the integral $\int_T^{2T}F_1^k(x)dx$\\
6.4 Upper bound of the integral $\int_T^{2T}F_1^{k-1}(x)F_2(x)dx$\\
6.5 Higher-power moments of $F_2(x)$\\
6.6 Evaluation of the integral $\int_T^{2T}F^k(x)dx$\\
6.7 Proof of Theorem 2\\
$\S 7$. Proofs of Theorem 3 and Theorem 4\\
7.1 Proof of Theorem 3 \\
7.2 Proof of  Theorem 4(the case $k=3$)\\
7.3 Proof of   Theorem 4(the case $k=4$)\\
$\S 8$. Proofs of Theorem 5 and Theorem 6\\

{\section{Introduction}}

Let $(M, g)$ be a closed $n$-dimensional Riemannian manifold with
metric $g$ and Laplace-Beltrami operator $\Delta.$ Let $N(t)$ denote
its spectral counting function , which is defined as the number of
the eigenvalues of $\Delta$ not exceeding $t.$ H\"ormander \cite{Ho}
proved that the Weyl's law
\begin{equation}
N(t)=\frac{vol(B_n)vol(M)}{(2\pi)^n}t^{n/2}+O(t^{(n-1)/2})
\end{equation}
holds, where $vol(B_n)$ is the volume of the $n$-dimensional unit
ball.

Let $$R(t)=N(t)-\frac{vol(B_n)vol(M)}{(2\pi)^n}t^{n/2}.$$
H\"ormander's estimate (1.1) in general is sharp , as the well-known
example of the sphere $S^n$ with its canonical metric shows
\cite{Ho}. However, it is a very difficult problem to determine the
optimal bound of $R(t)$ in any given manifold, which depends on the
properties of the associated geodesic flow. Many improvements have
been obtained for certain types of manifolds, see \cite{BG, BE,
BL,CPT, FR, G, Hu, IVr, KP, PT, Vo}.

\subsection{\bf The Weyl's law for ${\Bbb T}^2$: the Gauss circle problem}\

The simplest compact manifold with integrable geodesic flow is the
$2$-torus ${\Bbb T}^2={\Bbb R}^2/{\Bbb Z}^2.$ The exponential
functions $e(mx+ny)(m,n\in{\Bbb Z})$ form a basis of eigenfunctions
of the Laplace operator $\Delta=\partial_x^2+\partial_y^2,$ which
acts on functions on ${\Bbb T}^2.$ The corresponding eigenvalues are
$4\pi^2(m^2+n^2), m,n\in {\Bbb Z}.$ The spectral counting function
$$N_I(t)=\{\lambda_j\in Spec(\Delta): \lambda_j\leq t\}$$
is equal to the number of lattice points of ${\Bbb Z}^2$ inside a
circle of radius $\sqrt t/2\pi.$ The well-known Gauss circle problem
is to study the properties of the error term of the function
$N_I(t).$

In this case , the formula (1.1) becomes
\begin{equation}
N_I(t)=\frac{t}{4\pi}+O(t^{1/2}),
\end{equation}
which is the classical result of Gauss.  Let $R_I(t)$ denote the
error term in (1.2). Many authors improved the upper bound estimate
of $R_I(t).$ The latest result is due to Huxley\cite{Hu}, which
reads
\begin{eqnarray}
R_I(t)\ll t^{131/416}\log^{26957/8320} t.
\end{eqnarray}
Hardy\cite{Ha} conjectured that
\begin{equation}
R_I(t)\ll t^{1/4+\varepsilon}.
\end{equation}
Cram\'er\cite{Cr} proved that
\begin{eqnarray*}
\lim_{T\rightarrow \infty}T^{-3/2}\int_1^T|R_I(t)|^2dt=C, \ \
C=\frac{1}{6\pi^3}\sum_{n=1}^\infty\frac{r^2(n)}{n^{3/2}},
\end{eqnarray*}
which is consistent with Hardy's conjecture, where $r(n)$ denotes
the number of ways $n$ can be written as a sum of two squares.

Ivi\'c\cite{IVi} first used the large value technique to study the
higher power moments of $R_I(t).$ He proved that the estimate
\begin{eqnarray}
\int_1^T|R_I(t)|^{A}dt\ll T^{1+A/4+\varepsilon}
\end{eqnarray}
holds for each fixed $0\leq A\leq 35/8.$ The value of $A$ for which
(1.5) holds is closely related to the upper bound of $R_I(t).$ If we
insert the estimate (1.3) into Ivi\'c's machinery , we get that
(1.5) holds for $0\leq A\leq 262/27.$

Tsang\cite{Ts}  studied the third and the fourth moments of
$R_I(t).$ He proved that the following two asymptotic formulas are
true,
\begin{eqnarray}
&&\int_1^TR_I^3(t)dt=c_3T^{7/4} +O(T^{7/4- 1/14+\varepsilon}), \\&&
\int_1^TR_I^4(t)dt=c_4 T^{2}+O(T^{2-1/23+\varepsilon}),
\end{eqnarray}
where $c_3$ and $c_4$ are explicit constants.

Heath-Brown\cite{He} proved  that the function $t^{-1/4}R_I(t)$ has
a distribution function $f(\alpha)$   in the sense that , for any
interval $I$ we have
$$T^{-1} mes \{t\in [1,T]: t^{-1/4}R_I(t)\in I\}\rightarrow\int_If(\alpha)d\alpha$$
as $T\rightarrow \infty.$ He also proved  that for any real number
$k\in [0,9]$  the mean value
$$\lim_{T\rightarrow\infty}T^{-1-k/4}\int_1^T |R_I(t)|^kdt$$ converges to
a finite limit as $T$ tends to infinity. Moreover the same is true
for
$$\lim_{T\rightarrow\infty}T^{-1-k/4}\int_1^T R_I(t)^kdt$$
with $k=3,5,7,9.$

In \cite{Z1}, the author proved the following result: let $A>9$ be a
real number such that (1.5) holds , then for any integer $3\leq
k<A,$
 we have the asymptotic formula
\begin{equation}
\int_1^T R_I^k(t)dt= c_kT^{1+k/4}
+O(T^{1+k/4-\delta_k+\varepsilon}),
\end{equation}
where $c_k$ and $\delta_k>0$ are explicit constants. Especially, the
asymptotic formula (1.8) holds for $k=3,4,5,6,7,8,9.$

We remark that for $\delta_3$ we can take $\delta_3=7/20,$ which
  is due to Tsang for the third moment of the error
term of the Dirichlet divisor problem. However, Tsang didn't publish
this result.

For the fourth power moment of $R_I(t),$  the author\cite{Z2} proved
that we can take $\delta_4=3/28.$

\subsection{\bf The Weyl's law for $(2l+1)$-dimensional Heisenberg
manifold}\

Let $l\geq 1$ be a fixed integer and  $(H_l/\Gamma, g)$ be a
$(2l+1)$-dimensional Heisenberg manifold with a metric $g$. When
$l=1$, in \cite{PT} Petridis and Toth proved that
$R(t)=O(t^{5/6}\log t)$ for a special metric. Later in \cite{CPT}
this bound was improved to $O(t^{119/146+\varepsilon})$ for all
left-invariant Heisenberg metrics. For $l>1$   Khosravi and
Petridis\cite{KP} proved that $ R(t)=O(t^{l-7/41})$ holds for
rational Heisenberg manifolds. Both in \cite{CPT} and \cite{KP},
they first established a $\psi$-expression of $R(t)$ and then used
the van der Corput method of exponential sums  .
 Substituting Huxley's result of \cite{Hu} into the arguments of \cite{CPT} and \cite{KP}, we
 can get that the estimate
\begin{equation}
R(t)=O(t^{l-77/416}(\log t)^{26957/8320})
\end{equation}
holds for all rational $(2l+1)$-dimensional Heisenberg manifolds.

It was conjectured that for rational Heisenberg manifolds,  the
pointwise estimate
\begin{equation}
R(t)\ll t^{l-1/4+\varepsilon}
\end{equation}
holds, which was proposed  in  Petridis and Toth \cite{PT} for the
case $l=1$ and in Khosravi and Petridis\cite{KP} for the case $l>1.$
As an evidence of this conjecture, Petridis and Toth proved the
following $L^2$ result for $H_1$
\begin{eqnarray*}
\int_{I^3}\left|N(t;u)-\frac{1}{6\pi^2}vol(M(u))t^{3/2}\right|^2du\leq
C_\delta t^{3/2+\delta},
\end{eqnarray*}
where $I=[1-\varepsilon,1+\varepsilon].$ They also proved
$$\frac{1}{T}\int_T^{2T}\left|N(t)-\frac{1}{6\pi^2}vol(M)t^{3/2}\right|dt\gg T^{3/4}.$$

Now let $M=(H_l/\Gamma, g_l)$ be a $(2l+1)$-dimensional Heisenberg
manifold with the metric
 \[g_l:=\left(
\begin{array}{llcl}
 I_{2l\times 2l}&0\\
0&2\pi\\
\end{array}
\right),\] where $I_{2l\times 2l}$ is the identity matrix

M. Khosravi and John A. Toth\cite{KT} proved that
\begin{equation}
\int_1^T|R(t)|^2dt=C_{2,l}T^{2l+1/2}+O(T^{2l+1/4+\varepsilon}),
\end{equation}
where $C_{2,l}$ is an explicit constant .

 M. Khosravi \cite{K} proved that  the asymptotic formula
\begin{equation}
\int_1^TR^3(t)dt=C_{3,l}T^{3l+1/4}+O(T^{3l+3/14+\varepsilon})
\end{equation}
is true for some explicit constant $C_{3,l}.$

\bigskip

The aim of this  paper is to prove some power moment results for
$R(t),$ which are analogous to the results of $R_I(t)$ stated in
Section 1.1. The plan of this paper is as follows.  In Section 2 we
shall state our main results. In Section 3 we state some background
of the Heisenberg manifolds and give a new $ \psi$-expression of
$R(t).$ In Section 4 we quote some preliminary  Lemmas. We shall
prove our theorems in Section 5-7.

{\bf Notations.} For a real number  $t,$ let $[t]$ denote the
integer part of $t,$ $\{t\}=t-[t],$ $\psi(t)=\{t\}-1/2,$  $\Vert
t\Vert=\min(\{t\},1-\{t\}),$ $e(t)=e^{2\pi it}.$ $\varepsilon$
always denotes a sufficiently small positive constant. ${\Bbb C},
{\Bbb R},
 {\Bbb Z},  {\Bbb N}$ denote the set of complex numbers, the set of real numbers, the set of
integers, the set of positive integers, respectively. $n\sim M$
means that $N<n\leq 2N. $ $d(n)$ denotes the Dirichlet divisor
function , $r(n)$ denotes the number of ways $n$ can be written as a
sum of two squares, $\mu(n)$ denotes the M\"obius function.
Throughout this paper , ${\cal L}$ always denotes $\log T.$

\section{\bf Main results}

From now on, we always suppose that $R(t)$ denote the error term in
the Weyl's law for the $(2l+1)$-dimensional Heisenberg manifold
$(H_l/\Gamma, g_l).$

\subsection{\bf New results for power moments of $R(t)$}\

Our first result is the following Theorem 1, which is an analogue of
(1.5) for $R_I(t).$

{\bf Theorem 1.} Suppose $A\geq 0$ is a fixed real number. Then
\begin{eqnarray}
&&\int_1^T|R(t)|^Adt\ll T^{1+A(l-1/4)}{\cal L}^{41} \ \ \ (0\leq
A\leq 262/27),\\
&&\int_1^T|R(t)|^Adt\ll T^{A(l-1)+\frac{154+339A}{416}}{\cal
L}^{4A+1}\ \ \ (A>262/27).
\end{eqnarray}

{\bf Remark 2.1.} Let $A_0>1$ be a positive constant such that the
estimate
\begin{equation}
\int_1^T|R(t)|^{A_0}dt\ll T^{1+ A_0(l-1/4)+\varepsilon}
\end{equation}
holds, then Theorem 1 states that we can take $A_0=262/27.$ The
value of $A_0$ for which (2.3) holds is closely related to the upper
bound estimate of $R(t).$ The exponent $262/27$ follows from the
estimate (1.9). If the conjecture (1.10) were true, then obviously
(2.3) were true for any $A_0>0.$ Conversely, if the estimate (2.3)
were true for any $A_0>0,$ then we could show that the estimate
(1.10) were also true. The argument is as follows. The estimate
(2.3) implies that the upper bound
\begin{eqnarray*}
\int_{T/2}^T|R(2\pi x)|^{A_0}dx\ll T^{1+ A_0(l-1/4)+\varepsilon}.
\end{eqnarray*}
Suppose $|R(2\pi x)|$ reaches its largest value $V_0$ at some
$x_0\in [T/2,T].$ According to the formula (5.11), we have
\begin{eqnarray*}
V_0^{A_0+1}(x_0^{l-1/2}\log {x_0})^{-1}\ll \int_{T/2}^T|R(2\pi
x)|^{A_0}dx\ll T^{1+A_0(l-1/4) +\varepsilon},
\end{eqnarray*}
which implies that
$$V_0\ll T^{\frac{l+1/2+A_0(l-1/4)+2\varepsilon}{(A_0+1)}}.$$
Now the conjecture (1.10) follows from the above estimate by
choosing $A_0$ large.

\bigskip

Before stating our next theorem, we introduce some notations.

Suppose $f:{\Bbb N}\rightarrow {\Bbb R}$ is any function such that
$f(n)\ll n^\varepsilon,$ $k\geq 2$ is a fixed integer. Define
\begin{equation}
s_{k;v}(f):=\sum_{\sqrt{n_1}+\cdots
+\sqrt{n_v}=\sqrt{n_{v+1}}+\cdots +\sqrt{n_k }} \frac{f(n_1)\cdots
f(n_k)}{(n_1\cdots n_k)^{3/4}}\hspace{3mm}(1\leq v<k),
\end{equation}
\begin{equation}
B_k(f):=\sum_{v=1}^{k-1}{k-1\choose
v}s_{k;v}(f)\cos\frac{\pi(k-2v)}{4},
\end{equation}
\begin{equation}
\tau_l(n):=\sum_{ n=h(2r-h) }
\frac{e(lh/2)h^{1/2}}{(2r-h)^{1/2}}\left(1-\frac{h}{2r-h}\right)^{l-1}.
\end{equation}
We shall use $s_{k;v}(f)$  to denote both of the series (2.4) and
its value. The convergence of $s_{k;v}(f)$ was already proved in the
author\cite{Z1}. It is obvious that $|\tau_l(n)|\leq d(n).$

Suppose $A_0>2$ is a real number, define
\begin{eqnarray*}
&&K_0: =\min\{n\in {\Bbb N}:n\geq A_0, 2|n\},\\
&&s(k):=2^{k-2}+ (k-6)/4,\\
&&\sigma(k,A_0):=  \frac{A_0-k}{2(A_0-2)}, \\
&&\delta_1(k,A_0): =\sigma(k,A_0)/2s(K_0),\\
&&\delta_2(k,A_0): =\frac{\sigma(k,A_0)}{2s(k)+2\sigma(k,A_0)}.
\end{eqnarray*}

{\bf Theorem 2.} Let $A_0>9$ be a real number such that (2.3) holds
, then for any integer $3\leq k<A_0,$
 we have the asymptotic formula
\begin{equation}
\int_1^TR^k(t)dt=\frac{2^{3+5k/4-2kl}l^kB_k(\tau_l)}{(l!)^k\pi^{3k/4+kl}(4+k(4l-1))}T^{1+k(l-1/4)}
+O(T^{1+k(l-1/4)-\delta_1(k,A_0)+\varepsilon}).
\end{equation}

 {\bf Remark 2.2.} From  Theorem 1 we see that for $k\in
\{3,4,5,6,7,8,9\},$ we can get the asymptotic formula (2.7) with
$A_0=262/27$. Moreover,   the asymptotic (2.7) provides the exact
form of the main term for the integral $\int_1^TR^k(t)dt$ for $k\geq
10.$ From Theorem  2 we see that if the conjecture (1.10) were true,
then for all $k\geq 3$ we could get the asymptotic formula of
$\int_1^TR^k(t)dt.$

 When $3\leq k\leq 9,$ we have the following Theorem 3, which is
 better  than Theorem 2.

{\bf Theorem 3.}  For $3\leq k\leq 9,$ the asymptotic formula (2.7)
holds with $\delta_1(k,A_0)$ replaced by
 $\delta_2(k,262/27).$

{\bf Remark 2.3.} Note that $\delta_2(3,262/27)=181/1402$ and
$\delta_2(4,262/27)=11/230,$ which are small when comparison to the
corresponding results of $R_I(t)$.  The following   Theorem 4
improves these two exponents.

{\bf Theorem 4.}  We have
\begin{eqnarray}
\int_1^TR^3(t)dt=\frac{2^{27/4-6l}l^3B_3(\tau_l)}{(l!)^3\pi^{3l+9/4}(1+12l)
}T^{3l+1/4}+O(T^{3l+\varepsilon})
\end{eqnarray}
and
\begin{eqnarray}
\int_1^TR^4(t)dt=
\frac{2^{4-8l}l^3B_4(\tau_l)}{(l!)^4\pi^{4l+3}}T^{4l}+O(T^{4l-3/28+\varepsilon}).
\end{eqnarray}

The following two theorems are analogous to Heath-Brown's results
for $R_I(t).$

{\bf Theorem 5.} The function $t^{-(l-1/4)}R(t)$ has a distribution
function $f(\alpha)$   in the sense that , for any interval
$I\subset {\Bbb R}$ we have
$$T^{-1} mes \{t\in [1,T]: t^{-(l-1/4)}R(t)\in I\}\rightarrow\int_If(\alpha)d\alpha$$
as $T\rightarrow \infty.$ The function $f(\alpha)$ and its
derivatives satisfy
$$\frac{d^k}{d\alpha^k}f(\alpha)\ll_{A,k}(1+|\alpha|)^{-A} $$
for $k=0,1,2,\cdots$ and  $f(\alpha)$ can be extended to an entire
function.

{\bf Theorem 6.}  For any real number $k\in [0,262/27)$  the mean
value
$$\lim_{T\rightarrow\infty}T^{-1-k(l-1/4)}\int_1^T |R(t)|^kdt$$ converges to
a finite limit as $T$ tends to infinity.

{\bf Remark 2.4.}  Our approach of this paper implies a different
proof of the formula (1.11) and the error term estimate
$O(T^{2l+1/4+\varepsilon})$ therein can be improved slightly to
$O(T^{2l+1/4}{\cal L}^3).$ However, it is difficult to improve the
exponent $2l+1/4.$ For the constant  $C_{2,l}$ contained in (1.11),
we have a new expression,
namely,$$C_{2,l}=\frac{2^{9/2-4l}l^2B_2(\tau_l)}{(l!)^2\pi^{2l+3/2}(4l+1)}
=\frac{2^{9/2-4l}
l^2}{(l!)^2\pi^{2l+3/2}(4l+1)}\sum_{n=1}^{\infty}\frac{\tau_l^2(n)}{n^{3/2}},$$
 which is analogous to the expression of $C$ in Cram\'er's result for
the mean square of $R_I(t).$

\subsection{\bf Idea of the proof}\

M. Khosvari and John A. Toth\cite{KT} introduced an extra  parameter
and then used the Poisson summation formula to write the
corresponding error term in a form which can be estimated  by the
method of the stationary phase. Finally they eliminated the
parameter and got the asymptotic formula (1.11). This approach
worked very well for the mean square of $R(t).$ M. Khosvari\cite{K}
used the same approach to the third moment of $R(t)$ and proved
(1.12). But it is difficult to use this approach to the $k$-th power
moment when $k\geq 4.$

In this paper we shall use a  different approach to prove our
theorems. We can establish  an analogy between $R_I(t)$ and $R(t)$
and then we  prove our theorems by using this analogy.

It is well-known that $R_I(t)$ has the following truncated Voronoi's
formula
\begin{equation}
R_I(t)=-\frac{1}{2^{1/2}\pi^{3/2}}\sum_{n\leq
N}r(n)n^{-3/4}t^{1/4}\cos(2\sqrt{nt}+\pi/4)
+O(t^{1/2+\varepsilon}N^{-1/2})
\end{equation}
for $1\leq N\ll t,$ which follows from Lemma 3 of M\"uller\cite{Mu}.
The formula (2.10) plays an essential role in the proofs of  all
results stated in Section 1.1.

 $R(t)$  is a sum involving the row-teeth-function $\psi(u)$ and
doesn't have such a direct and easily usable Voronoi's formula at
hand.  However, by using the finite expression of $\psi(u)$ (see
Lemma 4.2)and van der Corput's B-process(see Lemma 4.3) we can prove
an formula(see Proposition 6.1) analogous to (2.10). This formula ,
although is weak when comparison to (2.10), is enough to prove all
our results.

 Nowak\cite{No} and Kuhleitner, Nowak\cite{Ku} first studied mean
squares of error terms involving $\psi(u),$ where they used another
finite Fourier expansion of $\psi(u)$, which is a  much more precise
form than Lemma 4.1. However, in our paper, we only use Lemma 4.1 to
study the upper bound of $R(t)$. For the asymptotic results of
$R(t),$ we shall use Lema 4.2, which seems more convenient to use.
The method of this paper can be used to study the power moments of
other error terms involving the function $\psi(u).$

\section{\bf Background of Heisenberg manifolds and the $\psi$-expression of $R(t)$}

In this section, we first review  some background of the Heisenberg
manifolds. The reader can see \cite{Fo}, \cite{GW} , \cite{St} for
more details. Finally, we give an $\psi$-expression of $R(t)$.

\subsection{\bf Heisenberg manifolds}\

 Suppose  $x\in {\Bbb R}^l$ is a row vector  and $y\in {\Bbb R}^l$ is a column vector. Define
 \[\gamma(x,y,t)=\left(
\begin{array}{llcl}
 1&x&t\\
0&I_l&y\\
0&0&1\\
\end{array}
\right), \ \ X(x,y,t)=\left(
\begin{array}{llcl}
 0&x&t\\
0&0&y\\
0&0&0\\
\end{array}
\right).\] The $(2l+1)$-dimensional Heisenberg group $H_l$ is
defined by
$$H_l=\{\gamma(x,y,t):  x, y\in {\Bbb R}^l, t\in {\Bbb R}\},$$
its Lie algebra is
$$\mathfrak{H}_l=\{X(x,y,t):  x, y\in {\Bbb R}^l, t\in {\Bbb R}\}.$$
We say $\Gamma$ is uniform discrete subgroup of $H_l$ if
$H_l/\Gamma$ is compact. A $(2l+1)$-dimensional Heisenberg manifold
is a pair $(H_l/\Gamma, g)$ for which $\Gamma$ is a uniform discrete
subgroup of $H_l$   and $g$ is a left $H_l$-invariant metric.

For every  $r$-tuple $(r_1,r_2,\cdots,r_l)\in {\Bbb N}^l$ such that
$r_j|r_{j+1}\ (j=1,2,\cdots,l-1)$, let $r{\Bbb Z}^l$ denote the
$l$-tuples $x=(x_1,x_2,\cdots,x_l)$ with $x_j\in r_j{\Bbb Z}$.
Define
$$\Gamma_r=\{\gamma(x,y,t): x\in r{\Bbb Z}^l, y\in r{\Bbb Z}^l, t\in {\Bbb Z}\}.$$
It is clear that $\Gamma_r$ is a  uniform discrete subgroup of
$H_l$.  According to Theorem 2.4 of \cite{GW}, the subgroup
$\Gamma_r$ classifies all the uniform discrete subgroups of $H_l$ up
to automorphisms. Thus (see \cite{GW}, Corollary 2.5) given any
Riemannian Heisenberg manifold $M=(H_l/\Gamma, g)$, there exists a
unique $l$-tuple $r$ as before and a left-invariant metric ${\tilde
g}$ on $H_l$ such that $M$ is isometric to $(H_l/\Gamma, {\tilde
g}).$ So (see \cite{GW}, 2.6(b)) we can replace the metric $g$ by
$\phi^*g,$ where $\phi$ is an inner automorphism such that the
direct sum split of the Lie algebra $\mathfrak{H}_l={\Bbb
R}^{2l}\oplus\mathfrak{Z}$ is orthogonal . Here  $\mathfrak{Z}$ is
the center of the Lie algebra and
 $${\Bbb R}^{2l}=\left\{\left(
\begin{array}{llcl}
 0&x&0\\
0&0&y\\
0&0&0\\
\end{array}
\right): x,y\in {\Bbb R}^l\right\}.$$ With respect to this
orthogonal split of $H_l$ the metric $g$ has the form \[\left(
\begin{array}{llcl}
 h& 0\\
0&g_{2l+1}\\
 \end{array}
\right),\] where $h$ is a positive-definite $2l\times 2l$ matrix and
$g_{2l+1}>0$ is a real number.

The volume of the Heisenberg manifold is given by
$$vol(H_l/\Gamma,g)=|\Gamma_r|\sqrt{det(g)}$$
with $|\Gamma_r|=r_1r_2\cdots r_l$ for $r=(r_1,r_2,\cdots, r_l).$

\subsection{\bf The spectrum of Heisenberg manifolds}\

Let $\Sigma$ be the spectrum of the Laplacian on $M=(H_l/\Gamma,
g_l),$ where the eigenvalues are counted with multiplicities.
According to \cite{GW}(P. 258), $\Sigma$ can be divided into two
parts $\Sigma_1$ and $\Sigma_2,$ where $\Sigma_1$ is the spectrum of
$2l$-dimensional torus
 and $\Sigma_2$ contains all eigenvalues of the form
 $$2\pi m^2+2\pi m(2n_1+\cdots+2n_l+l), m\in {\Bbb N}, n_j\in {\Bbb N}\cup\{0\},$$
 each eigenvalue counted with the multiplicity  $2m^l.$

\subsection{\bf the $\psi$-expression of $R(t)$}\

In this section we shall prove the following $\psi$-expression of
$R(t).$

{\bf Lemma 3.1.} We have
\begin{eqnarray}
&&R(2\pi x)=-\frac{4}{2^l(l-1)!}\sum_{1\leq m \leq \sqrt x}m(x-m^2
)^{l-1}\psi\left(\frac{x}{2m}-\frac{m}{2}-\frac{l}{2}\right)
+O(x^{l-1/2}).
\end{eqnarray}

\bigskip

Let
\begin{equation}
N(t)=N_I(t)+N_{II}(t),
\end{equation}
where
\begin{eqnarray*}
&&N_I(t):=\#\{\lambda: \lambda\in \Sigma_1, \lambda\leq t\},\\
&&N_{II}(t):=\#\{\lambda: \lambda\in \Sigma_2, \lambda\leq t\}.
\end{eqnarray*}

For $N_I(t)$ by H\"ormander's theorem we have
\begin{equation}
N_I(t)=\frac{1}{l!2^l}\left(\frac{t}{2\pi}\right)^{l}+O(t^{l-1/2}).
\end{equation}

\subsubsection{\bf Proof of Lemma 3.1 for $l=1$}\

Suppose $l=1.$ It is easily seen that
\begin{equation}
N_{II}(t)=\sum_{\stackrel{m^2+m(2n+1)\leq t/2\pi}{m\geq 1, n\geq
0}}2m.
\end{equation}
Thus we get
\begin{eqnarray*}
N_{II}(2\pi x)&&=\sum_{\stackrel{m^2+m(2n+1)\leq x}{m\geq 1, n\geq
0}}2m=\sum_{\stackrel{m^2+m(2n+1)\leq x}{1\leq m\leq \sqrt x, n\geq
0}}2m+O(\sqrt x)\\
&&=\sum_{1\leq m\leq \sqrt x}2m\sum_{0\leq n\leq
\frac{x}{2m}-\frac{m}{2}-\frac{1}{2}}+O(\sqrt x)\\
&&=\sum_{1\leq m\leq \sqrt x}2m\left[\frac{x}{2m}-\frac{m}{2}+\frac{1}{2}\right]+O(\sqrt x)\\
&&=\sum_{1\leq m\leq \sqrt
x}2m\left(\frac{x}{2m}-\frac{m}{2}\right)-\sum_{1\leq m\leq \sqrt
x}2m\psi\left(\frac{x}{2m}-\frac{m}{2}+\frac{1}{2}\right)+O(\sqrt
x)\\&& =\frac{2}{3}x^{3/2}-\frac{x}{2}-\sum_{1\leq m\leq \sqrt
x}2m\psi\left(\frac{x}{2m}-\frac{m}{2}+\frac{1}{2}\right)+O(\sqrt
x),
\end{eqnarray*}
which combined with  (3.3) for $l=1$ proves Lemma 3.1 of the case
$l=1.$

\subsubsection{\bf Proof of Lemma 3.1 for $l>1$}\

Suppose now $l>1.$ In this case for $N_{II}(t)$ we have
\begin{eqnarray}
N_{II}(t)&&=\sum_{\stackrel{m^2+m(2n_1+\cdots+2n_l+l)\leq t/2\pi}{m>0, n_j\geq 0}}2m^l\\
&&=\sum_{\stackrel{m^2+m(2n+l)\leq t/2\pi}{m>0,n\geq
0}}2m^l\sum_{\stackrel{n=n_1+\cdots+n_l}{n_j\geq 0}}1\nonumber\\
&&=\sum_{\stackrel{m^2+m(2n+l)\leq t/2\pi}{m>0,n\geq
0}}2m^l  {n+l-1\choose l-1} \nonumber\\
&&=\frac{2}{(l-1)!}\sum_{\stackrel{m^2+m(2n+l)\leq t/2\pi}{m>0,n\geq
0}}m^ln^{l-1}\nonumber\\&&\hspace{10mm}+\frac{l}{(l-2)!}\sum_{\stackrel{m^2+m(2n+l)\leq
t/2\pi}{m>0,n\geq 0}}m^ln^{l-2}+O(t^{l-1/2}).\nonumber\\
&&=\Sigma_{1,l}(t)+\Sigma_{2,l}(t)+O(t^{l-1/2}),\nonumber
\end{eqnarray}
say.

In order to evaluate $\Sigma_{1,l}(t)$ and $\Sigma_{2,l}(t)$, we
need the following

{\bf Lemma 3.2.}(Euler-Maclaurin summation formula) Suppose that
$f\in C^1[a,y]$, $a\in {\Bbb Z},$ then
$$\sum_{a\leq n\leq y}f(n)=\int_a^y\left(f(t)+\psi(t)f^{\prime}(t)\right)dt
-\psi(y)f(y)+\frac{f(a)}{2}.$$

\bigskip

Suppose $d\geq 0$ is a fixed integer. By Lemma 3.2 we get
\begin{eqnarray}
 \sum_{0\leq n\leq y}n^d=\left\{\begin{array}{ll}
y-\psi(y)+1/2,&\mbox{if $d=0$, }\\
\frac{y^{d+1}}{d+1}-\psi(y)y^{d}+O(y^{d-1}) ,& \mbox{if $d\geq 1.$}
\end{array}\right.
\end{eqnarray}

By (3.6) we get
\begin{eqnarray}
\Sigma_{1,l}(2\pi x)&&=\frac{2}{l!}\sum_{1\leq m\leq\sqrt
x}m^l\left(\frac{x}{2m}-\frac{m}{2}-\frac{l}{2}\right)^l\\
&&-\frac{2}{(l-1)!}\sum_{1\leq m\leq\sqrt
x}m^l\left(\frac{x}{2m}-\frac{m}{2}-\frac{l}{2}
\right)^{l-1}\psi\left(\frac{x}{2m}-\frac{m}{2}-\frac{l}{2}\right)\nonumber\\
&&+O(x^{l-1/2})\nonumber
\end{eqnarray}
and
\begin{eqnarray}
\Sigma_{2,l}(2\pi x)&&=\frac{l}{(l-1)!}\sum_{1\leq m\leq\sqrt
x}m^l\left(\frac{x}{2m}-\frac{m}{2}-\frac{l}{2}\right)^{l-1}
+O(x^{l-1/2}) .
\end{eqnarray}

Write
\begin{eqnarray}
\sum_{1\leq m\leq\sqrt
x}m^l\left(\frac{x}{2m}-\frac{m}{2}-\frac{l}{2}\right)^l
=2^{-l}\sum_{0\leq m\leq \sqrt x}(x-m^2-ml)^l-2^{-l}x^l.
\end{eqnarray}

Let $u(t)=(x-t^2-tl)^l.$ It is easy to check that
\begin{eqnarray*}
&&\int_0^{\sqrt
x}u(t)dt=\frac{(l!)^22^{2l}}{(2l+1)!}x^{l+1/2}-\frac{l}{2}x^l+O(x^{l-1/2}),\\
&&\int_0^{\sqrt x}\psi(t)u^{\prime}(t)dt\ll x^{l-1/2}.
\end{eqnarray*}

 By Lemma 3.2 we have
\begin{eqnarray}
 &&\ \ \ \ \ \ \ \sum_{0\leq m\leq \sqrt x}(x-m^2-ml)^l\\&&
 =\int_0^{\sqrt
 x}\left(u(t)+\psi(t)u^{\prime}(t)\right)dt-\psi(\sqrt x)u(\sqrt
 x)+\frac{u(0)}{2}\nonumber
\\ &&=\frac{(l!)^22^{2l}}{(2l+1)!}x^{l+1/2}-\frac{l-1}{2}x^l+O(x^{l-1/2}).\nonumber
\end{eqnarray}

From (3.9) and (3.10) we get
\begin{eqnarray}
\sum_{1\leq m\leq\sqrt
x}m^l\left(\frac{x}{2m}-\frac{m}{2}-\frac{l}{2}\right)^l=\frac{(l!)^22^{l}}{(2l+1)!}x^{l+1/2}-\frac{l+1}{2^{l+1}}x^l+O(x^{l-1/2}).
\end{eqnarray}

Write
\begin{eqnarray}
\sum_{1\leq m\leq\sqrt
x}m^l\left(\frac{x}{2m}-\frac{m}{2}-\frac{l}{2}\right)^{l-1}
=2^{-l+1}\sum_{1\leq m\leq \sqrt x}m(x-m^2-ml)^{l-1}.
\end{eqnarray}

Let $v(t)=t(x-t^2-tl)^{l-1}.$ It is easy to check that
\begin{eqnarray*}
&&\int_1^{\sqrt
x}v(t)dt= \frac{1}{2l}x^l+O(x^{l-1/2}),\\
&&\int_1^{\sqrt x}\psi(t)v^{\prime}(t)dt\ll x^{l-1}.
\end{eqnarray*}

 By Lemma 3.2 we have
\begin{eqnarray}
 &&\ \ \ \ \ \ \ \sum_{1\leq m\leq \sqrt x}m(x-m^2-ml)^{l-1}\\&&
 =\int_1^{\sqrt
 x}\left(v(t)+\psi(t)v^{\prime}(t)\right)dt-\psi(\sqrt x)v(\sqrt
 x)+\frac{v(1)}{2}\nonumber
\\ &&= \frac{1}{2l}x^l+O(x^{l-1/2}).\nonumber
\end{eqnarray}

From (3.12) and (3.13) we get
\begin{eqnarray}
\sum_{1\leq m\leq\sqrt
x}m^l\left(\frac{x}{2m}-\frac{m}{2}-\frac{l}{2}\right)^{l-1}=
\frac{1}{2^{l}l}x^l+O(x^{l-1/2}).
\end{eqnarray}

Since
$$\left(\frac{x}{2m}-\frac{m}{2}-\frac{l}{2}\right)^{l-1}-\left(\frac{x}{2m}-\frac{m}{2}\right)^{l-1}
\ll \left(\frac{x}{m} \right)^{l-2}$$ holds for $m\ll \sqrt x$,  we
have

\begin{eqnarray}
&&\ \ \ \ \ \ \ \sum_{m^2+ml\leq
x}m^l\left(\frac{x}{2m}-\frac{m}{2}-\frac{l}{2}
\right)^{l-1}\psi\left(\frac{x}{2m}-\frac{m}{2}-\frac{l}{2}\right)\\
&&=\sum_{m^2+ml\leq x}m^l\left(\frac{x}{2m}-\frac{m}{2}
\right)^{l-1}\psi\left(\frac{x}{2m}-\frac{m}{2}-\frac{l}{2}\right)+O(x^{l-1/2}).\nonumber\\
&&=\sum_{1\leq m \leq \sqrt{x}}m^l\left(\frac{x}{2m}-\frac{m}{2}
\right)^{l-1}\psi\left(\frac{x}{2m}-\frac{m}{2}-\frac{l}{2}\right)+O(x^{l-1/2}).\nonumber
\end{eqnarray}

Now Lemma 3.1 follows from (3.3), (3.5), (3.7), (3.8), (3.11) ,
(3.14) and (3.15) .

\subsection{\bf A weighted lattice point problem}\

For any positive integer $d,$ define
$$f_R(d):=\sum_{\stackrel{d=m(m+2n+l)}{m>0,n\geq 0}}m^ln^{l-1}.$$ From the proof of Lemma 3.1
it is easy to see that
\begin{equation}
R(2\pi x)=\frac{2}{(l-1)!}\sum_{d\leq
x}f_R(d)-\frac{2^{l+1}l!}{(2l+1)!}x^{l+1/2}-\frac{l+1}{2^ll!}x^{l}+O(x^{l-1/2}).
\end{equation}
So the evaluation of the counting function $N(2\pi x)$,  is
equivalent to  study the asymptotic behavior of the mean value
$\sum_{d\leq x}f_R(d)$.

\section{Some preliminary Lemmas}

We need the following Lemmas. Lemma 4.1 is due to Vaaler\cite{Va}.
Lemma 4.2 is well-known; see for example, Heath-Brown\cite{He2}.
Lemma 4.3 is Theorem 2.2 of Min\cite{Mi}, see also Lemma 6 of
Chapter 1 in \cite{V}. A weaker version of Lemma 4.3 can be found in
\cite{Ku}, which also suffices for our proof. Lemma 4.4 is Lemma 3.1
of the author\cite{Z1}. Lemma 4.5 is the first derivative test.
Lemma 4.6 is the famous Hal\'asz-Montgomery inequality, see for
example,  Ivi\'c\cite{IVi}.

{\bf Lemma 4.1.} Let $H\geq 2$ be any real number. Then
$$\psi(u)=\sum_{1\leq |h|\leq H}a(h)e(hu)+O(\sum_{0\leq |h|\leq H}b(h)e(hu)),$$
where $a(h)$ and $b(h)$ are functions such that $a(h)\ll 1/|h|,
b(h)\ll 1/H.$

{\bf Lemma 4.2.} Let $H\geq 2$ be any real number. Then
$$\psi(u)=-\sum_{1\leq |h|\leq H}\frac{e(hu)}{2\pi ih}+O\left(\min(1,\frac{1}{H\Vert u\Vert}) \right).$$

{\bf Lemma 4.3.} Suppose  $A_1,\cdots, A_5$ are absolute positive
constants, $f(x)$ and $ g(x)$ are algebraic functions in $[a,b]$ and
\begin{eqnarray*}
&&\frac{A_1}{R}\leq |f^{''}(x)|\leq\frac{A_2}{R},\  \ \ |f^{'''}(x)|\leq\frac{A_3}{RU},\ \ U\geq 1, \\
&&|g(x)|\leq A_4G,\ \ \  |g^{'}(x)|\leq A_5GU_1^{-1},\ \ U_1\geq 1,
\end{eqnarray*}
$[\alpha,\beta]$ is the image of  $[a,b]$ under the mapping
$y=f^{'}(x)$, then
\begin{eqnarray*}
\sum_{a<n\leq b}g(n)e(f(n))&=&e^{\pi i/4}\sum_{\alpha<u\leq
\beta}b_u
\frac{g(n_u)}{\sqrt{f^{''}(n_u)}}e\left(f(n_u)-un_u \right)\\
&&+O\left(G\log(\beta-\alpha+2)+G(b-a+R)(U^{-1}+U_1^{-1})\right)\\
&&+ O\left(G\min\left[\sqrt R,
\max\left(\frac{1}{<\alpha>},\frac{1}{<\beta>}\right)\right]\right),
\end{eqnarray*}
where  $n_u$  is the solution of $f^{'}(n)=u$,
\begin{eqnarray*}
<t>=\left\{\begin{array}{ll}
\Vert t\Vert,&\mbox{if $t$ not an integer,}\\
\beta-\alpha,& \mbox{if $t$ an integer,}
\end{array}\right.
\end{eqnarray*}
\begin{eqnarray*}
b_u=\left\{\begin{array}{ll}
1 ,&\mbox{if $\alpha<u<\beta$, or $\alpha, \beta$ not integers ,}\\
1/2,& \mbox{if $\alpha$ or $\beta$ are integers,}\\
\end{array}\right.
\end{eqnarray*}
\begin{eqnarray*}
\sqrt{f^{\prime\prime}}=\left\{\begin{array}{ll}
 \sqrt{f^{\prime\prime}},&\mbox{if $ f^{\prime\prime}>0,$}\\
 i\sqrt{|f^{\prime\prime}|},& \mbox{if $ f^{\prime\prime}<0. $}
\end{array}\right.
\end{eqnarray*}

{\bf Lemma 4.4.} Suppose $f(n)$ is an arithmetic function such that
$f(n)\ll n^\varepsilon$ , $1\leq v<k$ are fixed integers, $y>1$ is a
large parameter, and define
$$s_{k;v}(f;y): =
\sum_{\stackrel{\sqrt{n_1}+\cdots +\sqrt{n_v}=\sqrt{n_{v+1}}+\cdots
+\sqrt{n_k}} {n_1,\cdots,n_k\leq y}} \frac{f(n_1)\cdots
f(n_k)}{(n_1\cdots n_k)^{3/4}},1\leq v<k.$$ Then
$$|s_{k;v}(f)-s_{k;v}(f;y)|\ll y^{-1/2+\varepsilon}, 1\leq v<k.$$

{\bf Lemma 4.5.} Suppose $A,B\in {\Bbb R}, A\not= 0,$ then
$$\int_T^{2T}\cos(A\sqrt t+B)dt\ll  T^{1/2}|A|^{-1}.$$

{\bf Lemma 4.6.} Let $\mathcal {S}$ be an inner-product vector space
over ${\Bbb C},$ $(a,b)$ denote the inner product in $\mathcal {S}$
and $\Vert a\Vert_0=\sqrt{(a,a)}.$ Suppose that
$\xi,\varphi_1,\cdots,\varphi_R$ are arbitrary vectors in $\mathcal
{S}$. Then
$$\sum_{l\leq R}|(\xi,\varphi_l)|^2\leq \Vert
\xi\Vert_0^2 \max_{l_1\leq R}\sum_{l_2\leq
R}|(\varphi_{l_1},\varphi_{l_2})|.$$

\section{\bf Proof of Theorem 1 }

In order to prove Theorem 1 , we shall prove a large value estimate
of $R(x)$. For simplicity and convenience , we consider the function
$R_*(x):=R(2\pi x).$

\subsection{\bf Large value estimate of $R_*(x)$}\

In this subsection, we shall prove the following

 {\bf Theorem 5.1.} Suppose $T\leq x_1<x_2<\cdots <x_M\leq 2T$
satisfies $|R_{*}(x_s)|\gg V T^{l-1/2}(s=1,2,\cdots,M)$ and
$|x_j-x_i|\geq V\gg T^{7/32}{\cal L}^4(i\not= j),$ then we have
\begin{equation}
M\ll TV^{-3}{\cal L}^9+T^{15/4}V^{-12}{\cal L}^{41}.
\end{equation}

{\bf Proof.} Suppose $V<T_0$ is a parameter to be determined later.
Let $I$ be any subinterval of $[T,2T]$ of length not exceeding $T_0$
and let $G=I\cap \{x_1,x_2,\cdots,x_M\}.$ Without loss of
generality, suppose $G=\{x_1,x_2,\cdots,x_{M_0}\}.$

 Now let $J=[\frac{{\cal L}/2+\log {\cal L}-\log V}{2\log 2}],$ then by
 Lema 3.1 we have
\begin{eqnarray*}
R_{*}(x)&&=-\frac{4}{2^l(l-1)!}\sum_{j=0}^J \sum_{m\sim
\frac{x^{1/2}}{2^{j+1}}}m(x-m^2)^{l-1}\psi\left(\frac{x}{2m}-\frac{m}{2}+\frac{1}{2}\right)+O(\frac{V
T^{l-1/2}}{{\cal L}})\\
&&=-\frac{4}{2^l(l-1)!}\sum_{j_1=0}^{l-1}(-1)^{j_1}{l-1\choose j_1}
\sum_{j=0}^J \sum_{m\sim \frac{x^{1/2}}{2^{j+1}}}
x^{l-1-j_1}m^{2j_1+1}\psi\left(\frac{x}{2m}-\frac{m}{2}+\frac{1}{2}\right)\\&&\
\ \ \ \ \ \ \ \ \ +O(\frac{V T^{l-1/2}}{{\cal L}}).
\end{eqnarray*}
By Cauchy's inequality we get
$$R_{*}^2(x)\ll {\cal L}\sum_{j_1=0}^{l-1}\sum_{j=0}^J \left|\sum_{m\sim
\frac{x^{1/2}}{2^{j+1}}}x^{l-1-j_1}m^{2j_1+1}\psi\left(\frac{x}{2m}-\frac{m}{2}+\frac{1}{2}\right)\right|^2+\frac{V^2T^{2l-1}}{{\cal
L}^2}.$$ Summing over the set $G$ we get
\begin{eqnarray}
&&\ \ \ \ \ \ \ \ M_0V^2T^{2l-1}\ll \sum_{s\leq
M_0}|R_*(x_s)|^2\\&&\ll {\cal
L}\sum_{j_1=0}^{l-1}\sum_{j=0}^J\sum_{s\leq M_0} \left|\sum_{m\sim
\frac{x_s^{1/2}}{2^{j+1}}}x^{l-1-j_1}m^{2j_1+1}\psi\left(\frac{x_s}{2m}-\frac{m}{2}+\frac{1}{2}\right)\right|^2\nonumber\\
&&\ll {\cal L}^2 \sum_{s\leq M_0} \left|\sum_{m\sim
\frac{x_s^{1/2}}{2^{j+1}}}x^{l-1-j_1}m^{2j_1+1}\psi\left(\frac{x_s}{2m}-\frac{m}{2}+\frac{1}{2}\right)\right|^2\nonumber
\end{eqnarray}
for some fixed pair $(j_1,j)$ with $0\leq j_1\leq l-1, 0\leq j\leq
J$. For this fixed pair $(j_1,j),$ let $N=T^{1/2}2^{-j}$ and
$H=N^{2j_1+2}V^{-1}T^{-1/2-j_1}{\cal L}^2$. By Lemma 4.1 we get
\begin{equation}
M_0V^2T^{2l-1}\ll {\cal L}^2 \sum_{s\leq M_0} \left|\sum_{1\leq
h\leq H}c(h)\sum_{m\sim
\frac{x_s^{1/2}}{2^{j+1}}}x_s^{l-1-j_1}m^{2j_1+1}
e\left(-h\left(\frac{x_s}{2m}-\frac{hm}{2}\right)\right)\right|^2,
\end{equation}
where $c(h)$ is some function such that $c(h)\ll 1/h.$

For any integer $h>0$, define
$$S(x;h,j_1,j):=\sum_{m\sim \frac{x^{1/2}}{2^{j+1}}}x^{l-1-j_1}m^{2j_1+1}
e\left(-h\left(\frac{x}{2m}-\frac{hm}{2}\right)\right).$$

Take $$f(m)=-h\left(\frac{x}{2m}-\frac{hm}{2}\right),\ \
g(m)=x^{l-1-j_1}m^{2j_1+1}$$ in Lemma 4.3 we get
\begin{eqnarray}
&&\ \ \ \ \ \ \ \  \ \ S(x;h;j_1,j)\\&&=\frac{e^{\pi
i/4}}{i}\sum_{h(2^{2j-1}+\frac{1}{2})<r\leq
h(2^{2j+1}+\frac{1}{2})}\frac{b_rx^{l-\frac{1}{4}}h^{j_1+\frac 34}}{(2r-h)^{j_1+\frac 54}}e(-\sqrt{xh(2r-h)})\nonumber\\
&&\hspace{25mm}+O(T^{l-\frac 12}{\cal L}).\nonumber
\end{eqnarray}

Insert (5.4) into (5.3) we get
\begin{eqnarray}
M_0V^2T^{2l-1}&&\ll {\cal L}^2 \sum_{s\leq M_0}
|D(x_s;H,j_1,j)|^2+M_0T^{2l-1}{\cal
L}^4\\
&&\ll {\cal L}^2 \sum_{s\leq M_0} |D(x_s;H,j_1,j)|^2,\nonumber
\end{eqnarray}
where
$$D(x;H,j_1,j)=\sum_{1\leq h\leq
H}c(h) \sum_{h(2^{2j-1}+\frac{1}{2})<r\leq
h(2^{2j+1}+\frac{1}{2})}b_r\frac{x^{l-1/4}h^{j_1+3/4}}{(2r-h)^{j_1+5/4}}e(-\sqrt{xh(2r-h)}).$$

Let
$$\gamma(n;H,j_1,j)=\sum_{\stackrel{n=h(2r-h),1\leq h\leq H}{h(2^{2j-1}+\frac{1}{2})<r\leq
h(2^{2j+1}+\frac{1}{2})}}\frac{b_rc(h)h^{j_1+3/2}}{(2r-h)^{j_1+1/2}}$$
and let $N_0=H^2(2^{2j+1}+1/2).$  Then it is easy to see that
$\gamma(n;H,j_1,j)\ll d(n)$ and $N_0\ll TV^{-2}{\cal L}^4.$  Thus we
have
\begin{eqnarray}
M_0V^2&&\ll {\cal L}^2T^{1/2}\sum_{s\leq M_0}\left|\sum_{n\leq
N_0}\frac{\gamma(n;H,j_1,j)}{n^{3/4}}e(-\sqrt{nx_s})\right|^2\\
&&\ll {\cal L}^2T^{1/2}\sum_{s\leq M_0}\left|\sum_{v}\sum_{n\sim
N_02^{-v-1}}\frac{\gamma(n;H,j_1,j)}{n^{3/4}}e(-\sqrt{nx_s})\right|^2\nonumber\\
&&\ll {\cal L}^3 T^{1/2}\sum_{v}\sum_{s\leq M_0}\left|\sum_{n\sim
N_02^{-v-1}}\frac{\gamma(n;H,j_1,j)}{n^{3/4}}e(-\sqrt{nx_s})\right|^2\nonumber\\
&&\ll {\cal L}^4 T^{1/2}\sum_{s\leq M_0}\left|\sum_{n\sim
N_02^{-v-1}}\frac{\gamma(n;H,j_1,j)}{n^{3/4}}e(-\sqrt{nx_s})\right|^2\nonumber
\end{eqnarray}
for some $0\leq v\ll {\cal L},$ where in the third "$\ll$" we used
the Cauchy's inequality again. Let $N_1=N_02^{-v-1},$ then $N_1\ll
TV^{-2}{\cal L}^4.$

The procedure below is the same as the approach of the proof of
Theorem 13.8 in Ivi\'c\cite{IVi}, so we give only an outline. Take
$\xi=\{\xi_n\}_{n=1}^{\infty}$ with
$\xi_n=\gamma(n;H,j_1,j)n^{-3/4}$ for $n\sim N_1$ and zero
otherwise, and take $\varphi_s=\{\varphi_{s,n}\}_{n=1}^{\infty}$
 with $\varphi_{s,n}=e(-\sqrt{nx_s})$ for $n\sim N_1$ and zero
 otherwise. Then
 \begin{eqnarray*}
&&(\xi,\varphi_s)=\sum_{n\sim
N_1}\frac{\gamma(n;H,j_1,j)}{n^{3/4}}e(-\sqrt{nx_s}),\\
&&(\varphi_{l_1},\varphi_{l_2})=\sum_{n\sim
N_1}e\left(\sqrt{n}(\sqrt{x_{s_1}}-\sqrt{x_{s_2}})\right),\\
&& \Vert \xi\Vert_0^2=\sum_{n\sim
N_1}\frac{|\gamma(n;H,j_1,j)|^2}{n^{3/2}}\ll N_1^{-3/2}\sum_{n\sim
N_1}d^2(n)\ll N_1^{-1/2}{\cal L}^3.
 \end{eqnarray*}

By Lemma 4.6 we get
\begin{eqnarray}
&&M_0V^2\ll {\cal L}^7T^{1/2}N_1^{-1/2}\max_{s_1\leq
M_0}\sum_{s_2\leq M_0}\left|\sum_{n\sim
N_1}e\left(\sqrt{n}(\sqrt{x_{s_1}}-\sqrt{x_{s_2}})\right)\right|\\
&&\hspace{5mm}\ll {\cal L}^7T^{1/2}N_1^{1/2}+{\cal
L}^7T^{1/2}N_1^{-1/2}\max_{s_1\leq M_0}\sum_{s_2\leq M_0,s_2\not=
s_1}\left|\sum_{n\sim
N_1}e\left(\sqrt{n}(\sqrt{x_{s_1}}-\sqrt{x_{s_2}})\right)\right|.\nonumber
\end{eqnarray}

By the Kuzmin-Landau inequality and the exponent pair $(4/18,11/18)$
we get
\begin{eqnarray*}
\sum_{n\sim
N_1}e\left(\sqrt{n}(\sqrt{x_{s_1}}-\sqrt{x_{s_2}})\right)&&\ll
\frac{\sqrt
N_1}{|\sqrt{x_{s_1}}-\sqrt{x_{s_2}}|}+\left(\frac{|\sqrt{x_{s_1}}-\sqrt{x_{s_2}}|}{\sqrt{N_1}}\right)^{\frac{4}{18}}N_1^{\frac{11}{18}}\\
&&\ll
\frac{\sqrt{N_1T}}{|x_{s_1}-x_{s_2}|}+\left(\frac{|x_{s_1}-x_{s_2}|}{\sqrt{N_1T}}\right)^{\frac{4}{18}}N_1^{\frac{11}{18}}\\
&&\ll
\frac{\sqrt{N_1T}}{|x_{s_1}-x_{s_2}|}+T^{-1/9}T_0^{2/9}N_1^{1/2},
\end{eqnarray*}
where we used the mean value theorem and the estimate
$|x_{s_1}-x_{s_2}|\leq T_0.$

Insert this estimate into (5.7) we get
\begin{eqnarray}
M_0V^2&&\ll {\cal L}^7T^{1/2}N_1^{1/2}+{\cal
L}^7T^{1/2}N_1^{-1/2}\\
&&\hspace{15mm}\times\max_{s_1\leq M_0}\sum_{s_2\leq M_0,s_2\not=
s_1}\left(\frac{\sqrt{N_1T}}{|x_{s_1}-x_{s_2}|}+T^{-1/9}T_0^{2/9}N_1^{1/2}\right)\nonumber\\
&&\ll {\cal L}^7T^{1/2}N_1^{1/2}+{\cal L}^7TV^{-1}+{\cal
L}^7M_0T^{1/2-1/9}T_0^{2/9}\nonumber\\
&&\ll {\cal L}^9TV^{-1}+{\cal L}^7M_0T^{7/18}T_0^{2/9},\nonumber
\end{eqnarray}
where we used the facts that $\{x_s\}$ is $V$-spaced and $N_1\ll
TV^{-2}{\cal L}^4.$ Take $T_0=V^9T^{-7/4}{\cal L}^{-32},$  it is
easy to check that $T_0\gg V$ if $V\gg T^{7/32}{\cal L}^4.$ We get
for this $T_0$ that
$$M_0\ll {\cal L}^9TV^{-3}.$$

Now we divide the interval  $[T,2T]$ into  $O(1+T/T_0)$ subintervals
of length not exceeding $T_0.$ In each interval of this type, the
number of $x_s$ is  at most $O({\cal L}^9TV^{-3})$ . So we have
\begin{equation}
M\ll {\cal L}^9TV^{-3}(1+\frac{T}{T_0})\ll {\cal L}^9TV^{-3}+{\cal
L}^{41}T^{15/4}V^{-12}.
\end{equation}

This completes the proof of Theorem 5.1.

\subsection{\bf Proof of Theorem 1}\

Now we prove Theorem 1. When $A=0,$ Theorem 1 is trivial . When
$0<A<2,$ it follows from (1.11) and the H\"older's inequality. So
later we always suppose $A>2$. It suffices for us to prove the
estimate
\begin{eqnarray}
\int_T^{2T}|R_*(x)|^Adx\ll\left\{\begin{array}{ll}
T^{1+\frac{3A}{4}+A(l-1)}\log^{41} T,&\mbox{if $2<A\leq 262/27,$}\\
T^{\frac{154+339A}{416}+A(l-1)}\log^{4A+1} T,& \mbox{if $A>262/27.$}
\end{array}\right.
\end{eqnarray}

Suppose $x^\varepsilon<y\leq x/2,$  by (3.16) we get
\begin{eqnarray*}
&&\ \ \ \ \ \ \ \ \ \ |R_*(x+y)-R_*(x)|\\&&\ll \sum_{x<n\leq
x+y} f_R(n)+|(x+y)^{l+1/2}-x^{l+1/2}|+|(x+y)^{l}-x^{l}|\\
&&\ll x^{l-1/2}\sum_{x<n\leq x+y}d(n)+x^{l-1/2}y
\\&& \ll x^{l-1/2}y\log x,
\end{eqnarray*}
where we used the well-known estimate
$$\sum_{x<n\leq x+y}d(n)\ll y\log x$$
and the obvious bounds $  f_R(n)\ll n^{l-1/2}d(n).$ So there exists
an absolute constant $c_0$ such that
\begin{eqnarray*}
|R_*(x+y)-R_*(x)|\leq c_0 x^{l-1/2}y\log x,
\end{eqnarray*}
which impies that  if $|R_*(x)|\geq 2c_0 x^{l-1/2}y\log x,$ then
\begin{eqnarray}
|R_*(x+y)|\geq |R_*(x)|-|R_*(x+y)-R_*(x)|\geq c_0 x^{l-1/2}y\log x.
\end{eqnarray}

From (5.11) and a similar argument  to (13.70) of Ivi\'c\cite{IVi}
we may write
\begin{equation}
\int_T^{2T}|R_*(x)|^Adx\ll T^{1+A(l-\frac 14) }{\cal L}
+\sum_{V}V\sum_{r\leq N_V}|R_*(x_r)|^A,
\end{equation}
where $T^{1/4}\leq V=2^m\leq T^{131/416}{\cal L}^4,$ $V
T^{l-1/2}<|R_*(x_r)|\leq 2VT^{l-1/2}$
 $(r=1,\cdots, N_V)$ and $|t_r-t_s|\geq V$ for $r\not= s\leq N=N_V.$

If $2<A\leq 11,$ then by Theorem 5.1   we get (recall $V\ll
T^{131/416}{\cal L}^4$)
\begin{eqnarray}
&&V\sum_{r\leq N_V}|R_*(x_r)|^A\ll T^{A(l-1/2) }N_VV^{A+1}
\\&&\ll T^{A(l-1/2)}\left({\cal L}^9TV^{A-2}+{\cal
L}^{41}T^{\frac{15}{4} }V^{A-11}\right)\nonumber\\
&&\ll T^{A(l-1/2)}\left(T^{1+\frac{131}{416}(A-2)}{\cal L}^{4A+1}+
T^{1+\frac{A}{4} }{\cal L}^{41}\right)\nonumber\\
&&\ll T^{A(l-1/2)}\left(T^{\frac{154+131A}{416} }{\cal L}^{4A+1}+
T^{1+\frac{A}{4} }{\cal L}^{41}\right).\nonumber
\end{eqnarray}

If $A> 11,$ then by Theorem 5.1    we get
\begin{eqnarray}
&&V\sum_{r\leq N_V}|R_*(x_r)|^A\ll T^{A(l-1/2)}N_VV^{A+1}
\\&&\ll T^{A(l-1/2)}\left({\cal L}^9T V^{A-2}+{\cal
L}^{41}T^{\frac{15}{4} }V^{A-11}\right)\nonumber\\
&&\ll T^{A(l-1/2)}\left(T^{1 +\frac{131}{416}(A-2)}{\cal L}^{4A+1}+
T^{\frac{15}{4} +\frac{131}{416}(A-11) }{\cal L}^{4A-3}\right)\nonumber\\
&&\ll T^{A(l-1/2)+\frac{154+131A}{416} }{\cal L}^{4A+1} .\nonumber
\end{eqnarray}
Now (5.10) follows from (5.12)-(5.14) by noting that
$(154+131A)/416\leq 1+A/4$ for $2<A\leq 262/27$ and $(154+131A)/416>
1+A/4$ for $A> 262/27$.

\section{\bf Proof of Theorem 2 }

In this section we shall prove Theorem 2. Suppose that $3\leq
k<A_0,$ $A_0>9$ is a fixed  real number such that (2.3) holds. It
suffices for us to evaluate the integral $\int_T^{2T}R_*^k(x)dx$ ,
where $T$ is a large real number.

Suppose $H$ is a large parameter to be determined later. By Lemma
3.1 and Lemma  4.2 we have
\begin{eqnarray}
&&R_*(x)=F(x)+O(T^{l-1/2}G(x)),\\
&&F(x)=\frac{2^{1-l}}{(l-1)!\pi i}\sum_{1\leq |h|\leq
H}\frac{1}{h}\sum_{m\leq
\sqrt x}m(x-m^2)^{l-1}e\left(h\left(\frac{x}{2m}-\frac{m}{2}+\frac{l}{2}\right)\right),\nonumber\\
&&G(x)= \sum_{m\leq \sqrt{2T}}\min\left(1,\frac{1}{H\Vert
\frac{x}{2m}-\frac{m}{2}+\frac{l}{2}\Vert}\right).\nonumber
\end{eqnarray}

\subsection{\bf The upper bound of  $\int_T^{2T}G(x)dx$}\

In this subsection we shall prove the following Lemma 6.1.

{\bf Lemma 6.1.} We have
\begin{equation}
\int_T^{2T}G(x)dx\ll T^{3/2}H^{-1}{\cal L}.
\end{equation}

{\bf Proof.} Obviously we have
\begin{eqnarray}
G(x)&&\ll  \sum_{2m\leq \sqrt{2T}}\min\left(1,\frac{1}{H\Vert
\frac{x}{4m}+\frac{1}{2}\Vert}\right)\\
&&\hspace{5mm}+ \sum_{2m-1\leq \sqrt{2T}}\min\left(1,\frac{1}{H\Vert
\frac{x}{2(2m-1)}\Vert}\right)\nonumber\\
&&\ll  G_1(x)+G_2(x),\nonumber
\end{eqnarray}
where
\begin{eqnarray*}
G_1(x)&&=\sum_{m\leq 4\sqrt{T}}\min\left(1,\frac{1}{H\Vert
\frac{x}{m}+\frac{1}{2}\Vert}\right),\\
G_2(x)&&=\sum_{m\leq 4\sqrt{T}}\min\left(1,\frac{1}{H\Vert
\frac{x}{m}\Vert}\right).
\end{eqnarray*}
Thus
\begin{equation}
\int_T^{2T}G(x)dx\ll \int_T^{2T}G_1(x)dx+ \int_T^{2T}G_2(x)dx.
\end{equation}

 We have
\begin{eqnarray}
\int_T^{2T}G_1(x)dx&& \ll \sum_{m\leq 4\sqrt T}
\int_T^{2T}\min\left(1,\frac{1}{H \Vert
\frac{x}{m}+\frac{1}{2}\Vert }\right)dx\\
&&\ll  \sum_{m\leq 4\sqrt
T}m\int_{\frac{T}{m}}^{\frac{2T}{m}}\min\left(1,\frac{1}{H \Vert
x+\frac{1}{2}\Vert }\right)dx\nonumber\\
&& \ll  T\sum_{m\leq 4\sqrt T} \int_0^{1/2}\left(1,\frac{1}{H \Vert
x+\frac{1}{2}\Vert }\right)dx\nonumber\\
&&\ll T^{3/2}\int_0^{1/2}\left(1,\frac{1}{H \Vert x+\frac{1}{2}\Vert
}\right)dx\nonumber \\
&&\ll T^{3/2}H^{-1}{\cal L}.\nonumber
\end{eqnarray}

Similarly we have
\begin{equation}
\int_T^{2T}G_2(x)dx\ll   T^{3/2}H^{-1}{\cal L}.
\end{equation}
Now Lemma 6.1 follows from (6.3)-(6.6).

\subsection{\bf An analogue of Voronoi's formula for $R_*(x)$}\

In this subsection we shall give an analogue of  Voronoi's formula
for $R_*(x)$. Our main tool is Lemma 4.2 and Lemma 4.3.

We begin with $F(x).$
 From the definition of $F(x)$ we have
 \begin{eqnarray}
F(x)&&=\frac{2^{1-l}}{(l-1)!\pi i}\sum_{1\leq |h|\leq H}\frac
1h\sum_{m\leq
\sqrt{x}}m(x-m^2)^{l-1}e\left(h\left(\frac{x}{2m}-\frac
m2-\frac l2\right)\right)\\
&&=\frac{2^{1-l}}{(l-1)!\pi
i}\sum_{j_1=0}^{l-1}(-1)^{j_1}{l-1\choose j_1}\sum_{1\leq |h|\leq
H}\frac 1h\sum_{m\leq \sqrt{x}}x^{l-1-j_1}m^{2j_1+1}\nonumber\\
&&\ \ \ \ \ \ \ \ \ \ \ \ \ \ \ \ \ \ \ \ \ \ \ \ \ \ \  \times
e\left(h\left(\frac{x}{2m}-\frac m2-\frac
l2\right)\right)\nonumber\\
&&=\frac{2^{1-l}}{(l-1)!}\sum_{j_1=0}^{l-1}(-1)^{j_1}{l-1\choose
j_1}F(x;j_1),\nonumber
 \end{eqnarray}
say, where
$$F(x;j_1):=\frac{1}{\pi
i} \sum_{1\leq |h|\leq H}\frac 1h\sum_{m\leq
\sqrt{x}}x^{l-1-j_1}m^{2j_1+1} e\left(h\left(\frac{x}{2m}-\frac
m2-\frac l2\right)\right).$$

Let $J=[({\cal L}-\log {\cal L})/2\log 2]$   we get
\begin{eqnarray}
&&\ \ \ \ \ \ \ \ \ \ \ F(x;j_1)\\&&=\frac{1}{\pi i}\sum_{-H\leq
h\leq -1}\frac{ e(-lh/2)}{h}\sum_{j=0}^J\sum_{m\sim \sqrt
x2^{-j-1}}x^{l-1-j_1}m^{2j_1+1}e\left(h\left(\frac{x}{2m}-\frac{m}{2}\right)\right)\nonumber \\
&&\hspace{4mm}+\frac{1}{\pi i}\sum_{1\leq h\leq
H}\frac{e(-lh/2)}{h}\sum_{j=0}^J\sum_{m\sim \sqrt
x2^{-j-1}}x^{l-1-j_1}m^{2j_1+1}e\left(h\left(\frac{x}{2m}-\frac{m}{2}\right)\right)+O(x^{l-1}{\cal
L}^2)\nonumber\\
&&= -\frac{1}{\pi i}\sum_{1\leq h\leq H}\frac{
e(lh/2)}{h}\sum_{j=0}^J\sum_{m\sim \sqrt
x2^{-j-1}}x^{l-1-j_1}m^{2j_1+1}e\left(-h\left(\frac{x}{2m}-\frac{m}{2}\right)\right)\nonumber\\
&&\hspace{4mm} +\frac{1}{\pi i}\sum_{1\leq h\leq
H}\frac{e(-lh/2)}{h}\sum_{j=0}^J\sum_{m\sim \sqrt
x2^{-j-1}}x^{l-1-j_1}m^{2j_1+1}e\left(h\left(\frac{x}{2m}-\frac{m}{2}\right)\right)+O(x^{l-1}{\cal
L}^2) \nonumber\\
&&=-\frac{\Sigma_3}{\pi i}+\frac{\overline{\Sigma_3}}{\pi i}
+O(x^{l-1}{\cal L}^2),\nonumber
\end{eqnarray}
where
\begin{eqnarray*}
\Sigma_3=\sum_{1\leq h\leq H}\frac{ e(lh/2)}{h}\sum_{j=0}^J
S(x;h,j_1,j)
\end{eqnarray*}
with $S(x;h,j_1,j)$  defined in Section 5.1.

Insert the formula (5.4) into $\Sigma_3$ we get
\begin{eqnarray}
&& \Sigma_3=e^{-\pi i/4}\sum_{1\leq h\leq
H}\frac{e(lh/2)}{h}\sum_{j=0}^J\sum_{h(2^{2j-1}+\frac{1}{2})<r\leq
h(2^{2j+1}+\frac{1}{2})}\\
&&\hspace{20mm}\frac{b_rx^{l-1/4}h^{j_1+3/4}}{(2r-h)^{j_1+5/4}}e(-\sqrt{xh(2r-h)})+O(T^{l-1/2}{\cal L}^3)\nonumber\\
&&=e^{-\pi i/4}\sum_{1\leq h\leq
H}\frac{e(lh/2)}{h}\sideset{}{^{\prime}}  \sum_{h<r\leq
h(2^{2J+1}+1/2)}\frac{x^{l-1/4}h^{j_1+3/4}}{(2r-h)^{j_1+5/4}}e(-\sqrt{xh(2r-h)})\nonumber\\
&&\hspace{15mm}+O( T^{l-1/2}{\cal L}^3)\nonumber,
\end{eqnarray}
where in the last sum the symbol "$\prime$" means that if $h$ is an
even integer , then the term $r=h(2^{2J+1}+1/2)$ should be halved .

Inserting (6.9) into (6.8) we get
\begin{eqnarray}
&& \ \ \ \ \ \ F(x;j_1)=\frac{1}{\pi i}\sum_{1\leq h\leq
H}\frac{e(lh/2)}{h}\sideset{}{^{\prime}}  \sum_{h<r\leq
h(2^{2J+1}+1/2)}\frac{x^{l-1/4}h^{j_1+3/4}}{(2r-h)^{j_1+5/4}}\\
&&\ \ \ \ \ \ \ \ \ \ \ \ \ \ \  \times
\left(e(\sqrt{xh(2r-h)}+\frac
18)-e(\sqrt{xh(2r-h)}-\frac 18)\right)+O(x^{l-1/2}{\cal L}^3)\nonumber\\
&&=\frac{2x^{l-1/4}}{\pi }\sum_{1\leq h\leq
H}\frac{e(lh/2)}{h}\sideset{}{^{\prime}}  \sum_{h<r\leq
h(2^{2J+1}+1/2)}\frac{
h^{j_1+3/4}}{(2r-h)^{j_1+5/4}}\sin(2\pi\sqrt{xh(2r-h)}+\frac{\pi}{4})\nonumber\\
&&\ \ \ \ \ \ \ \ \ \ \ \ \ \ \ +O(x^{l-1/2}{\cal L}^3)\nonumber\\
&&=\frac{2x^{l-1/4}}{\pi }\sum_{1\leq h\leq
H}\frac{e(lh/2)}{h}\sideset{}{^{\prime}}  \sum_{h<r\leq
h(2^{2J+1}+1/2)}\frac{
h^{j_1+3/4}}{(2r-h)^{j_1+5/4}}\cos(2\pi\sqrt{xh(2r-h)}-\frac{\pi}{4})\nonumber\\
&&\ \ \ \ \ \ \ \ \ \ \ \ \ \ \ +O(x^{l-1/2}{\cal L}^3).\nonumber
\end{eqnarray}

From (6.7) and (6.10) we get
\begin{eqnarray}
  F(x)&&=\frac{2^{2-l}x^{l-1/4}}{(l-1)!\pi}\sum_{1\leq h\leq H}
\sideset{}{^{\prime}}  \sum_{h<r\leq
h(2^{2J+1}+1/2)}\cos(2\pi\sqrt{xh(2r-h)}-\frac{\pi}{4})\nonumber\\
&& \ \ \ \ \ \ \  \ \ \ \ \ \ \ \ \ \ \times
\sum_{j_1=0}^{l-1}(-1)^{j_1}{l-1\choose j_1}\frac{
e(lh/2)h^{j_1-1/4}}{(2r-h)^{j_1+5/4}}+O(x^{l-1/2}{\cal
L}^3)\nonumber\\
&&=\frac{2^{2-l}x^{l-1/4}}{(l-1)!\pi}\sum_{1\leq h\leq H}
\sideset{}{^{\prime}}  \sum_{h<r\leq
h(2^{2J+1}+1/2)}\frac{e(lh/2)}{h^{1/4}(2r-h)^{5/4}}\left(1-\frac{h}{2r-h}\right)^{l-1}\nonumber\\
&&\ \ \ \ \ \ \ \ \ \ \ \ \times
\cos(2\pi\sqrt{xh(2r-h)}-\frac{\pi}{4})+O(x^{l-1/2}{\cal
L}^3)\nonumber
\end{eqnarray}

Define
\begin{eqnarray*}
\tau_l(n;H,T):&&=\sum_{\stackrel{n=h(2r-h),1\leq h\leq H}{h<r\leq
h(2^{2J+1}+1/2)}}\frac{e(lh/2)h^{1/2}}{(2r-h)^{1/2}}\left(1-\frac{h}{2r-h}\right)^{l-1}.
\end{eqnarray*}
 We then have
\begin{equation}
F(x)=\frac{2^{2-l}x^{l-1/4}}{(l-1)!\pi}\sum_{1\leq n\leq
H^2(2^{2J+1}+1/2)}\frac{\tau_l(n;H,T)}{n^{3/4}}\cos\left(2\pi\sqrt{xn}-\frac{\pi}{4}\right)+O(x^{l-1/2}{\cal
L}^3),
\end{equation}
where the "$\prime$" terms  for which $2|h$ and $ r=h(2^{2J+1}+1/2)$
are absorbed into the error term.

Obviously that $|\tau_l(n;H,T)|\leq d(n).$ By the definition of $J$
we see that if $n\leq \min(H, T{\cal L}^{-2}),$ then
$\tau_l(n;H,T)=\tau_l(n),$ where $\tau_l(n)$ was defined in Section
2.1.

From (6.1) and (6.11) we get the following Proposition, which is an
analogue of the Voronoi's formula.

{\bf Proposition 6.1.} Suppose $T\leq x\leq 2T, H\geq T, J=[({\cal
L}-\log {\cal L})/2\log 2].$ Then we have
\begin{eqnarray}
R_*(x)&&=\frac{2^{2-l}x^{l-1/4}}{(l-1)!\pi}\sum_{1\leq n\leq
H^2(2^{2J+1}+1/2)}\frac{\tau(n;H,T)}{n^{3/4}}\cos\left(2\pi\sqrt{xn}-\frac{\pi}{4}\right)\\
&&\hspace{20mm}+O(T^{l-1/2}G(x)+T^{l-1/2}{\cal L}^3),\nonumber
\end{eqnarray}
where $\tau_l(n;H,T)=\tau_l(n)$ for $n\leq T{\cal L}^{-2}.$

\bigskip

Now suppose $T^{\varepsilon}<y\leq T{\cal L}^{-2}$ is a parameter to
be determined. Define
\begin{eqnarray}
&&F_1(x):= \theta_lx^{l-1/4} \sum_{1\leq n\leq
 y}\frac{\tau_l(n)}{n^{3/4}}\cos\left(2\pi\sqrt{xn}-\frac{\pi}{4}\right),\\
 &&F_2(x):=F(x)-F_1(x),\nonumber
\end{eqnarray}
where $\theta_l:=\frac{2^{2-l}}{(l-1)!\pi}.$

\subsection{\bf Evaluation of the integral $\int_T^{2T}F_1^k(x)dx$}\

The $k$-th power moment of $F_1(x)$ provides the main term in
Theorem 2, so in this subsection we shall evaluate the integral
$\int_T^{2T}F_1^k(x)dx$.  For simplicity we set ${\Bbb I}=\{0,1\}$
and  $${\Bbb N}^k=\{{\bf n}: {\bf n}=(n_1,\cdots, n_k),n_j\in {\Bbb
N},1\leq j\leq k\}.$$ For each element ${\bf i}=(i_1,\cdots,
i_{k-1})\in {\Bbb I}^{k-1},$ put $|{\bf i}|=i_1+\cdots+i_{k-1}.$ By
the elementary formula
\begin{equation}
\cos a_1\cdots \cos a_k= \frac{1}{2^{k-1}}\sum_{{\bf i}\in {\Bbb
I}^{k-1}} \cos (a_1+(-1)^{i_1}a_2+(-1)^{i_2}a_3+\cdots
+(-1)^{i_{k-1}}a_k)\end{equation} we have
\begin{eqnarray*}
F_1^k(x)&&=\theta_l^{k}x^{ k(l-\frac 14)} \sum_{n_1\leq y}\cdots
\sum_{n_k\leq y}\frac{\tau_l(n_1)\cdots \tau_l(n_k)}{(n_1\cdots n_k
)^{3/4}}
\prod_{j=1}^k\cos(2\pi\sqrt{n_jx}-\frac{\pi}{4})\\
&&= \frac{\theta_l^k x^{k(l-1/4)}}{2^{k-1}} \sum_{{\bf i}\in {\Bbb
I}^{k-1} } \sum_{n_1\leq y}\cdots \sum_{n_k\leq y}
\frac{\tau_l(n_1)\cdots \tau_l(n_k)}{(n_1\cdots n_k)^{3/4}} \cos
(2\pi\sqrt x\alpha({\bf n};{\bf i})-\frac{\pi}{4} \beta({\bf i})),
\end{eqnarray*}
where
\begin{eqnarray*}
\alpha({\bf n};{\bf i}): &&=
\sqrt{n_1}+(-1)^{i_1}\sqrt{n_2}+(-1)^{i_2}\sqrt{n_3}+\cdots
+(-1)^{i_{k-1}}\sqrt
{n_k},\\
\beta({\bf i}): &&=1+(-1)^{i_1}+(-1)^{i_2}+\cdots +(-1)^{i_{k-1}}.
\end{eqnarray*}
Thus we can write
\begin{equation}
F_1^k(x)=\frac{\theta_l^k}{2^{k-1}}(S_1(x)+S_2(x)),
\end{equation}
where
\begin{eqnarray*}
S_1(x):&&= x^{k(l-1/4)}\sum_{{\bf i}\in {\Bbb I}^{k-1}
}\cos\left(-\frac{\pi\beta({\bf i}) }{4}\right)
\sum_{\stackrel{n_j\leq y,1\leq j\leq k}{\alpha({\bf n};{\bf i})=
0}}
\frac{\tau_l(n_1)\cdots \tau_l(n_k)}{(n_1\cdots n_k)^{3/4}},\\
S_2(x):&&= x^{k(l-1/4)}\sum_{{\bf i}\in {\Bbb I}^{k-1} }
\sum_{\stackrel{n_j\leq y,1\leq j\leq k}{\alpha({\bf n};{\bf
i})\not= 0}} \frac{\tau_l(n_1)\cdots \tau_l(n_k)}{(n_1\cdots
n_k)^{3/4}} \cos\left(2\pi\alpha({\bf n};{\bf i})\sqrt
x-\frac{\pi\beta({\bf i})}{4}\right) .\end{eqnarray*}

First consider the contribution of $S_1(x).$ We have
\begin{equation}
\int_T^{2T}S_1(x)dx= \sum_{ {\bf i}\in {\Bbb
I}^{k-1}}\cos\left(-\frac{\pi\beta({\bf i})}{4}\right)
\sum_{\stackrel{n_j\leq y,1\leq j\leq k}{\alpha({\bf n};{\bf i})=
0}} \frac{\tau_l(n_1)\cdots \tau_l(n_k)}{(n_1\cdots
n_k)^{3/4}}\int_T^{2T}x^{k(l-1/4) }dx.
\end{equation}

 It is easily seen that if $\alpha({\bf n};{\bf i})=0,$ then $1\leq |{\bf i}|\leq k-1.$
So
\begin{equation}
\sum_{\stackrel{n_j\leq y,1\leq j\leq k}{\alpha({\bf n};{\bf i})=
0}} \frac{\tau_l(n_1)\cdots \tau_l(n_k)}{(n_1\cdots
n_k)^{3/4}}=s_{k;|{\bf i}|}(\tau_l;y).\end{equation}

By Lemma 4.4 we get
\begin{equation}
\int_T^{2T}S_1(x)dx=B_k^{*}(\tau_l)\int_T^{2T}x^{k(l-1/4)}dx
+O(T^{1+k(l-1/4)+\varepsilon}y^{-1/2}),
\end{equation}
where
$$B_k^{*}(\tau_l):=
\sum_{{\bf i}\in {\Bbb I}^{k-1} }\cos\left(-\frac{\pi\beta({\bf
i})}{4}\right) \sum_{\stackrel{ {\bf n}\in {\Bbb N}^k}{\alpha({\bf
n};{\bf i})= 0}} \frac{\tau_l(n_1)\cdots \tau_l(n_k)}{(n_1\cdots
n_k)^{3/4}}.$$

For any ${\bf i}\in {\Bbb I}^{k-1}\setminus {\bf 0},$ let
\begin{eqnarray*}
S(\tau_l;{\bf i}): &&=\sum_{\stackrel{{\bf n}\in {\Bbb
N}^k}{\alpha({\bf n};{\bf i})= 0}} \frac{\tau_l(n_1)\cdots \tau(n_k)
}{(n_1\cdots n_k)^{3/4}} .
\end{eqnarray*}
It is easily seen that if $ |{\bf i}|= |{\bf i}^{\prime}|$ or $|{\bf
i}|+|{\bf i}^{\prime}| =k,$ then
$$S(\tau_l;{\bf i})=S(\tau_l;{\bf i}^{\prime})=s_{k;|{\bf i}|}(\tau_l).$$ From
$(-1)^j=1-2j(j=0,1)$ we also have $\beta({\bf i} )=k-2|{\bf i}| .$
So we get
\begin{eqnarray}
B_k^{*}(\tau_l)&&=\sum_{v=1}^{k-1}\sum_{|{\bf i}| =v}
\cos(-\frac{\pi\beta({\bf i})}{4})S(\tau_l;{\bf i})\\&&
=\sum_{v=1}^{k-1}s_{k;v}(\tau_l)\cos \frac{\pi(k-2v)}{4}\sum_{|{\bf
i}|
=v}1\nonumber\\
&&=\sum_{v=1}^{k-1}{k-1\choose v}s_{k;v}(\tau_l)\cos
\frac{\pi(k-2v)}{4}=B_k(\tau_l).\nonumber
\end{eqnarray}

Now we consider the contribution of $S_2(x).$ By Lemma 4.5 we get
\begin{eqnarray}
&&\int_T^{2T}S_2(x)dx\ll T^{1/2+k(l-1/4)}U_k(y),\\
&&U_k(y):=
 \sum_{{\bf i}\in {\Bbb I}^{k-1}
} \sum_{\stackrel{n_j\leq y,1\leq j\leq k}{\alpha({\bf n};{\bf
i})\not= 0}} \frac{d(n_1)\cdots d(n_k)}{(n_1\cdots
n_k)^{3/4}|\alpha({\bf n};{\bf i})|}.\nonumber
\end{eqnarray}

In \cite{Z1} the author proved
\begin{equation}
U_k(y)\ll y^{s(k) +\varepsilon},
\end{equation}
 where $s(k)$ was defined in Section 2.1.  Thus
\begin{equation}
\int_T^{2T}S_2(x)dx\ll T^{1/2+k(l-1/4)+\varepsilon}y^{s(k)}.
\end{equation}

Hence from (6.15)-(6.22) we get

{\bf Lemma 6.2.} For fixed $k\geq 3,$ we have
\begin{eqnarray}
\int_T^{2T}F_1^k(x)dx&&=\frac{2^{1+k-lk}l^kB_k(\tau_l)}{(l!)^k\pi^k}\int_T^{2T}
x^{k(l-\frac 14)}dx \\
&&\ \ \ \ \ +O(T^{1+k(l-\frac 14)
+\varepsilon}y^{-1/2}+T^{1/2+k(l-\frac
14)+\varepsilon}y^{s(k)}).\nonumber
\end{eqnarray}

\subsection{\bf Upper bound of the integral $\int_T^{2T}F_1^{k-1}(x)F_2(x)dx$}\

In this subsection we shall estimate the integral
$\int_T^{2T}F_1^{k-1}(x)F_2(x)dx$.

Let $$N_2:=H^2(2^{2J+1}+1/2), J:=[({\cal L}-\log {\cal L})/2\log
2].$$
  By (6.14)  we have
\begin{eqnarray*}
F_1^{k-1}(x)F_2(x)&&=\theta_l^{k}x^{k(l-1/4) } \sum_{y<n_1\leq
N_2}\sum_{n_2\leq y}\cdots \sum_{n_k\leq
y}\frac{\tau_l(n_1;H,T)\tau_l(n_2)\cdots \tau_l(n_k)}{(n_1\cdots n_k
)^{3/4}}\\&&\hspace{25mm}\times
\prod_{j=1}^k\cos(2\pi\sqrt{n_jx}-\frac{\pi}{4})\\
&&=\frac{\theta_l^{k}x^{k(l-1/4) }  }{2^{k-1}} \sum_{{\bf i}\in
{\Bbb I}^{k-1} } \sum_{y<n_1\leq N_2}\sum_{n_2\leq y}\cdots
\sum_{n_k\leq y}
\frac{\tau_l(n_1;H,T)\tau_l(n_2)\cdots \tau_l(n_k)}{(n_1\cdots n_k)^{3/4}}\\
&&\hspace{30mm}\times \cos (2\pi\sqrt x\alpha({\bf n};{\bf i}
)-\frac{\pi}{4} \beta( {\bf i})).
\end{eqnarray*}
  Thus
\begin{equation}
F_1^{k-1}(x)F_2(x)= \frac{\theta_l^{k}  }{2^{k-1}}(S_3(x)+S_4(x)),
\end{equation}
where
\begin{eqnarray*}
&&S_3(x)= x^{k(l-1/4)}\sum_{{\bf i}\in {\Bbb I}^{k-1}
}\cos\left(-\frac{\pi\beta({\bf i})}{4}\right)\sum_{y<n_1\leq N_2}
\sum_{\stackrel{n_j\leq y,2\leq j\leq k}{\alpha({\bf n};{\bf i})=
0}}
\frac{\tau_l(n_1;H,T)\tau_l(n_2)\cdots \tau_l(n_k)}{(n_1\cdots n_k)^{3/4}},\\
&&S_4(x)= x^{k(l-1/4)}\sum_{{\bf i}\in {\Bbb I}^{k-1}
}\sum_{y<n_1\leq N_2} \sum_{\stackrel{n_j\leq y,2\leq j\leq
k}{\alpha({\bf n};{\bf i})\not= 0}}
\frac{\tau_l(n_1;H,T)\tau_l(n_2)\cdots \tau_l(n_k)}{(n_1\cdots
n_k)^{3/4}}
\\&&\hspace{47mm}\times\cos\left(2\pi\alpha({\bf n};{\bf i})\sqrt x-\frac{\pi\beta({\bf i})}{4}\right) .
\end{eqnarray*}

By Lemma 4.4 we have
\begin{eqnarray}
&&  \ \ \ \ \ \ \ \ \ \int_T^{2T}S_3(x)dx\\&&\ll\sum_{ {\bf i}\in
{\Bbb I}^{k-1}} \sum_{y<n_1\leq N_2}\sum_{\stackrel{n_j\leq y,2\leq
j\leq k}{\alpha({\bf n};{\bf i})= 0}} \frac{d(n_1)d(n_2)\cdots
d(n_k)}{(n_1\cdots n_k)^{3/4}}\int_T^{2T}x^{k(l-1/4) }dx\nonumber\\
&&\ll T^{1+k(l-1/4)}\sum_{ {\bf i}\in {\Bbb I}^{k-1}}
\sum_{y<n_1\leq N_2}\sum_{\stackrel{n_j\leq y,2\leq j\leq
k}{\alpha({\bf n};{\bf i})= 0}} \frac{d(n_1)d(n_2)\cdots
d(n_k)}{(n_1\cdots
n_k)^{3/4}}\nonumber\\
&&\ll T^{1+k(l-1/4)}\sum_{v=1}^{k-1}|s_{k;v}(d;y)-s_{k;v}(d)|\ll
T^{1+k(l-1/4)+\varepsilon}y^{-1/2}.\nonumber
\end{eqnarray}

Now we consider the contribution of $S_4(x).$ By Lemma 4.5 we get
\begin{eqnarray}
&&\int_T^{2T}S_4(x)dx\ll T^{1/2+k(l-1/4)}(\Sigma_4+\Sigma_5),\\
&&\Sigma_4= \sum_{ {\bf i}\in {\Bbb I}^{k-1} } \sum_{y<n_1\leq
k^2y}\sum_{\stackrel{n_j\leq y,2\leq j\leq k}{\alpha({\bf n};{\bf
i})\not= 0}}
\frac{d(n_1)\cdots d(n_k)}{(n_1\cdots n_k)^{3/4}|\alpha({\bf n};{\bf i})|},\nonumber\\
&&\Sigma_5=  \sum_{{\bf i}\in {\Bbb I}^{k-1} } \sum_{k^2y<n_1\leq
N_2}\sum_{\stackrel{n_j\leq y,2\leq j\leq k}{\alpha({\bf n};{\bf
i})\not= 0}} \frac{d(n_1)\cdots d(n_k)}{(n_1\cdots
n_k)^{3/4}|\alpha({\bf n};{\bf i})|}.\nonumber
\end{eqnarray}

By (6.21) we have
\begin{eqnarray}
\Sigma_4\ll U_k(k^2y)\ll y^{s(k)+\varepsilon}.
\end{eqnarray}

When $n_1>k^2y,$ it is easy to show that $|\alpha({\bf n};{\bf
i})|\gg n_1^{1/2},$ which implies that
\begin{eqnarray}
\Sigma_5&&\ll \sum_{k^2y<n_1\leq N_2}\sum_{ n_j\leq y,2\leq j\leq k}
\frac{d(n_1)\cdots d(n_k)}{(n_2\cdots n_k)^{3/4}n_1^{5/4}} \ll
y^{\frac{k-2}{4}}{\cal L}^k
\end{eqnarray}
if noting that
$$\sum_{n\leq y}d(n)n^{-3/4}\ll y^{1/4}\log y,\
\ \sum_{n> y}d(n)n^{-5/4}\ll y^{-1/4}\log y.$$

  From (6.24)-(6.28) we have

\begin{equation}
\int_T^{2T}F_1^{k-1}(x)F_2(x)dx \ll
T^{1+k(l-1/4)+\varepsilon}y^{-1/2}+T^{1/2+k(l-1/4)+\varepsilon}y^{s(k)}
.
\end{equation}

\subsection{\bf Higher-power moments of $F_2(x)$}\

In this subsection we shall study the higher power moments of
$F_2(x).$ From now on, we take $H: =T^{A_0}.$

We first study the mean-square of $ F_2(x).$ Recall that
$N_2=H^2(2^{2J+1}+1/2).$ Since $\tau(n,H;T)=\tau(n)$ for $n\leq
y\leq T{\cal L}^{-2},$ we have
 \begin{eqnarray*}
F_2(x) &&= \theta_l x^{l-1/4} \sum_{y<n\leq N_2}\frac{
\tau_l(n;H,T)}{n^{3/4}} \cos(2\pi\sqrt{nx}-\pi/4)+O(T^{l-1/2}{\cal
L}^3)
 \\
&&\ll T^{l-1/4}\left|\sum_{y<n\leq N_2}\frac{\tau_l(n;H,T)
}{n^{3/4}} e(2\sqrt{nx})\right| +T^{l-1/2}{\cal L}^3 ,
\end{eqnarray*}
which implies
\begin{equation}
\int_T^{2T}F_2^2(x)dx \ll  T^{2l-1/2}\int_T^{2T}\left| \sum_{y<n\leq
N_2}\frac{ \tau(n;H,T)}{n^{ 3/4}}
e(2\sqrt{nx})\right|^2dx+T^{2l}{\cal L}^3
\end{equation}
\begin{eqnarray*}
&&\ll  T^{2l+1/2}\sum_{y<n\leq N_2}\frac{d^2(n)}{n^{3/2}}
+T^{2l}\sum_{y<m<n\leq N_2}\frac{d(n)d(m)}{(mn)^{3/4}(\sqrt n-\sqrt m)}\\
&&\ll  \frac{T^{2l+1/2} {\cal L}^3}{y^{1/2}}  ,
\end{eqnarray*}
where we used the estimate
 $\sum_{n\leq u}d^2(n)\ll u\log^3 u$ and the well-known Hilbert's
 inequality.

Now suppose $y$ satisfies $y^{2b(K_0)}\leq T.$ Hence from Lemma 6.2
we get that
$$ \int_T^{2T}| F_1(x)|^{K_0}dx\ll T^{1+K_0(l-1/4)+\varepsilon},$$
which implies
\begin{equation}
 \int_T^{2T}| F_1(x)|^{A_0}dx\ll T^{1+A_0(l-1/4)+\varepsilon}
\end{equation}
since $A_0\leq K_0.$ Trivially we have $$G(x)= \sum_{m\leq
\sqrt{2T}}\min\left(1,\frac{1}{H\Vert
\frac{x}{2m}-\frac{m}{2}+\frac{1}{2}\Vert}\right)\ll T^{1/2}.$$ By
Lemma 6.1 we have
\begin{equation}
\int_T^{2T}G^{A_0}(x)dx\ll T^{(A_0-1)/2}\int_T^{2T}G(x)dx\ll
T^{A_0/2+1}H^{-1}{\cal L}\ll T{\cal L}.
\end{equation}

From (2.3), (6.31) and (6.32) we get
\begin{equation}
 \int_T^{2T}| F_2(x)|^{A_0}dx\ll
 \int_T^{2T}(|R_*(x)|^{A_0}+| F_1(x)|^{A_0}+T^{A_0(l-1/2)}G^{A_0}(x))dx\ll T^{1+A_0(l-1/4)+\varepsilon}.
\end{equation}

For any $2<A<A_0,$ from (6.30), (6.33) and  H\"older's inequality we
get that
\begin{equation}
 \int_T^{2T}| F_2(x)|^{A}dx=
\int_T^{2T}|F_2(x)|^{\frac{2(A_0-A)}{A_0-2}+\frac{A_0(A-2)}{A_0-2}}dx
\end{equation}
$$\ll \left( \int_T^{2T}F_2^2{x}dx\right)^{\frac{A_0-A}{A_0-2}}
 \left( \int_T^{2T}|F_2(x)|^{A_0}dx\right)^{\frac{A-2}{A_0-2}}
\ll T^{1+A(l-\frac 14) +\varepsilon}y^{-\frac{A_0-A}{2(A_0-2)}}.$$

Namely, we have the following Lemma 6.3.

{\bf Lemma 6.3.} Suppose $T^\varepsilon\leq y\leq T^{1/2b(K_0)},$
$2<A<A_0,$  then
\begin{equation}
 \int_T^{2T}|F_2(x)|^{A}dx\ll
T^{1+A(l-\frac{1}{4})+\varepsilon}y^{-\frac{A_0-A}{2(A_0-2)}}.
\end{equation}

\subsection{\bf  Evaluation of the integral $\int_T^{2T}F^k(x)dx$.}\

Suppose $3\leq k<  A_0$ and  $T^\varepsilon\leq y\leq
T^{1/2b(K_0)}.$ By the elementary formula
$$(a+b)^k=a^k+ka^{k-1}b+O(|a^{k-2}b^2|+|b|^k)$$ we get
\begin{eqnarray}
\int_T^{2T}F^k(x)dx&&=\int_T^{2T}F_1^kdx
+k\int_T^{2T}F_1^{k-1}(x)F_2(x)dx\\
&&+O\left(\int_T^{2T}|F_1^{k-2}(x)F_2^2(x)|dx+\int_T^{2T}|F_2(x)|^kdx\right).\nonumber
\end{eqnarray}

By  (6.31), Lemma 6.3 and H\"older's inequality we get
\begin{eqnarray}
&&\hspace{15mm}\int_T^{2T}|F_1^{k-2}(x)F_2^2(x)|dx \\&&\ll
\left(\int_T^{2T}|F_1(x)|^{A_0}dx\right)^{\frac{k-2}{A_0}}
\left(\int_T^{2T}|F_2(x)|^{\frac{2A_0}{A_0-k+2}}dx\right)^
{\frac{A_0-k+2}{A_0}}\nonumber\\
&& \ll T^{1+ k(l-\frac
14)+\varepsilon}y^{-\frac{A_0-k}{2(A_0-2}}\nonumber.\end{eqnarray}

 Now from (6.29), (6.37), Lemma 6.2 and Lemma 6.3($A=k$) we get
\begin{eqnarray}
&&\ \ \ \ \ \ \ \ \ \
\int_T^{2T}F^k(x)dx\\&&=\frac{2^{1+k-lk}l^kB_k(\tau_l)}{(l!)^k\pi^k}
\int_T^{2T}x^{k(l-1/4)}d
x  +O(T^{1+k(l-1/4)- \frac{A_0-k}{4(A_0-2)s(K_0)}+\varepsilon} )\nonumber\\
&&=
\frac{2^{1+k-lk}l^kB_k(\tau_l)}{(l!)^k\pi^k}\int_T^{2T}x^{k(l-1/4)}d
x +O(T^{1+k(l-1/4)-
 \delta_1(k,A_0)+\varepsilon} )\nonumber
 \end{eqnarray}
by choosing $y=T^{1/2s(K_0)}.$

\subsection{\bf Proof of Theorem 2}\

By Proposition 6.1 and the elementary formula
$(a+b)^k=a^k+O(|a|^{k-1}|b|+|b|^k)$ we get
\begin{eqnarray}
R_*^k(x)&&=F^k(x)+O(|F(x)|^{k-1}T^{l-1/2}G(x)+|F(x)|^{k-1}
T^{l-1/2}{\cal
L}^3)\\&&\hspace{7mm}+O(T^{k(l-1/2)}G^k(x)+T^{k/2}{\cal
L}^{3k}).\nonumber
\end{eqnarray}

By Lemma 6.1 and H\"older's inequality we get
\begin{eqnarray}
&&\int_T^{2T}|F(x)|^{k-1}T^{l-1/2}G(x)dx
\\&&\ll T^{l-1/2}\left(\int_T^{2T}|F(x)|^{A_0}dx\right)^{\frac{k-1}{A_0}}
\left(\int_T^{2T}G(x)^{\frac{A_0}{A_0-k+1}}dx\right)^{\frac{A_0-k+1}{A_0}}\nonumber\\
&&\ll T^{l-1/2}\
T^{(1+A_0(l-\frac{1}{4})+\varepsilon)\frac{k-1}{A_0}}\left(T^{\frac{k-1}{2(A_0-k+1)}}\int_T^{2T}G(x)dx\right)^{\frac{A_0-k+1}{A_0}}\nonumber\\
&&\ll T^{l-1/2}\
T^{(1+\frac{3A_0}{4}+\varepsilon)\frac{k-1}{A_0}}\left(T^{\frac{k-1}{A_0-k+1}}
T^{-(A_0-2)}\right)^{\frac{A_0-k+1}{A_0}}\nonumber
\\&&\ll T^{1/4+k(l-1/4)}.\nonumber
\end{eqnarray}

By  H\"older's inequality  again we get
\begin{eqnarray}
&&\int_T^{2T}|F(x)|^{k-1}T^{l-1/2}{\cal L}^3dx
\\&&\ll T^{l-1/2}{\cal L}^3\left(\int_T^{2T}|F(x)|^{A_0}dx\right)^{\frac{k-1}{A_0}}
  T^{\frac{A_0-k+1}{A_0}}\nonumber\\
&&\ll  T^{1+k(l-1/4)-1/4+\varepsilon}.\nonumber
\end{eqnarray}

From (6.39)-(6.41) we get
\begin{eqnarray*}
\int_T^{2T}R_*^k(x)dx&&=
\frac{2^{1+k-lk}l^kB_k(\tau_l)}{(l!)^k\pi^k}\int_T^{2T} x^{k(l-\frac
14)}dx +O(T^{1+k(l-1/4)-
 \delta_1(k,A_0)+\varepsilon} ),
 \end{eqnarray*}
which implies that
\begin{eqnarray}
&&\ \ \ \ \ \ \ \int_1^TR_*^k(x)dx\\&&=
\frac{2^{1+k-lk}l^kB_k(\tau_l)}{(l!)^k\pi^k}\int_1^{T} x^{k(l-\frac
14)}dx +O(T^{1+k(l-1/4)-
 \delta_1(k,A_0)+\varepsilon} )\nonumber\\
&&=
\frac{2^{1+k-lk}l^kB_k(\tau_l)}{(l!)^k\pi^k(k(l-1/4)+1)}T^{1+k(l-1/4)}+O(T^{1+k(l-1/4)-
 \delta_1(k,A_0)+\varepsilon} ).\nonumber
 \end{eqnarray}

Now Theorem 2 follows from (6.42) if noting that
\begin{eqnarray}
&&  \int_1^TR^k(t)dt=2\pi\int_1^{\frac{T}{2\pi}}R_*^k(x)dx+O(1) .
\end{eqnarray}

\section{\bf Proofs of Theorem 3 and Theorem 4}

In this section we  prove Theorem 3 and Theorem 4. Throughout this
section , let $A_0=262/27$ and $H=T^{262/27}.$

\subsection{\bf Proof of Theorem 3}\

The proof of Theorem 3 is almost the same as that of Theorem 2. So
we give only an outline.

Using the argument of the proof of Theorem 5.1 to $F_1(x)$ directly,
we can show that the estimate
\begin{equation}
\int_T^{2T}|F_1(x)|^{A_0}dx\ll T^{1+A_0(l-1/4)}{\cal L}^{50}
\end{equation}
holds for $y\leq T^{77/208}.$ Here we remark that if we want to get
the result of the type (7.1) we have to assume $yT^{2\theta}\ll T$
when recalling the formula (5.8),  where $\theta=131/416.$

From (7.1), (6.32) and (2.3) we get
\begin{equation}
\int_T^{2T}|F_2(x)|^{A_0}dx\ll T^{1+A_0(l-1/4)+\varepsilon} .
\end{equation}

From (7.1) and (7.2) we can show that Lemma 6.2 holds for $y\leq
T^{77/208}.$ Using other estimates in the proof of Theorem 2, we can
get by choosing $y=T^{(k-2)/2(A_0-2)s(k)}$ that the estimate
\begin{eqnarray*}
\int_T^{2T}R_*^k(x)dx&&=
\frac{2^{1+k-lk}l^kB_k(\tau_l)}{(l!)^k\pi^k}\int_T^{2T} x^{k(l-\frac
14)}dx+O(T^{1+k(l-\frac{1}{4})+\varepsilon}y^{-\frac{A_0-k}{2(A_0-2)}}  )\\
&&\hspace{20mm}+O(T^{1/2+k(l-1/4)+\varepsilon}y^{s(k)})\\
&&= \frac{2^{1+k-lk}l^kB_k(\tau_l)}{(l!)^k\pi^k}\int_T^{2T}
x^{k(l-\frac 14)}dx
+O(T^{1+k(l-1/4)-\delta_2(k,262/27)+\varepsilon})
 \end{eqnarray*}
holds, which implies Theorem 3.

\subsection{\bf Proof of Theorem 4(the case $k=3$)}\

 By  the argument of Section 6.3 we get
\begin{eqnarray}
 \int_T^{2T}F_1^3(x)dx&&=\frac{2^{4-2l}l^3B_3(\tau_l)}{(l!)^3\pi^3}\int_T^{2T}x^{3(l-
\frac{1}{4})}dx\\&&+O(T^{3l+1/4 +\varepsilon}y^{-1/2})
+O(T^{3l-1/4}U_3(y)),\nonumber
\end{eqnarray}
where $U_3(y)$ was defined in Section 6.3. In Lemma 2.6 of
\cite{Z3}, the author proved
 \begin{equation}
U_3(y)\ll  y^{1/4+\varepsilon}.
\end{equation}

By the argument of Section 6.4 we get
\begin{eqnarray}
\int_T^{2T}F_1^2(x)F_2(x)dx \ll
T^{3l+1/4+\varepsilon}y^{-1/2}+T^{3l-1/4 }y^{\frac{1}{4}}{\cal L}^3
+T^{3l-1/4}U_3^*(y),
\end{eqnarray}
where
$$U_3^*(y)= \sum_{ {\bf i}\in {\Bbb I}^2 } \sum_{y<n_1\leq
9y}\sum_{\stackrel{n_2\leq y, n_3\leq y}{\alpha({\bf n};{\bf
i})\not= 0}} \frac{d(n_1)d(n_2)d(n_3)}{(n_1n_2n_3)^{3/4}|\alpha({\bf
n};{\bf i})|}.$$

By (7.4) we get
\begin{equation}
U_3^*{y}\ll U_3(9y)\ll y^{1/4+\varepsilon}.
\end{equation}

Trivially we have $F_1(x)\ll x^{l-1/4}y^{1/4}{\cal L}.$ So by (6.30)
we get that
\begin{equation}
\int_T^{2T}|F_1(x)F_2^2(x)|dx\ll T^{l-1/4}y^{1/4}{\cal
L}\int_T^{2T}| F_2^2(x)|dx\ll T^{3l+1/4}y^{-1/4}{\cal L}^4
\end{equation}
holds for $y\leq T .$   By the trivial estimate $F_2(x)\ll x^l$ and
(6.30) we get
\begin{eqnarray}
\int_T^{2T}|F_2^3(x)|dx\ll T^{l} \int_T^{2T}|F_2(x)|^2dx\ll
T^{3l+1/2}y^{-1/2}{\cal L}^3.
\end{eqnarray}

From (6.36) with $k=3$ and   (7.3)-(7.8) with $y=T{\cal L}^{-2}$ we
get
\begin{eqnarray*}
\int_T^{2T}F^3(x)dx=\frac{2^{4-2l}l^3B_3(\tau_l)}{(l!)^3\pi^3}
\int_T^{2T}x^{3(l-1/4)}dx +O(T^{3l+\varepsilon}),
\end{eqnarray*}
which combining the arguments of Section 6.7  gives
\begin{eqnarray}
\int_T^{2T}R_*^3(x)dx=\frac{2^{4-2l}l^3B_3(\tau_l)}{(l!)^3\pi^3}
\int_T^{2T}x^{3(l-1/4)}dx +O(T^{3l+\varepsilon}).
\end{eqnarray}
The formula (2.8) follows from (7.9).

\subsection{\bf Proof of Theorem 4(the case $k=4$)}\

 By  the argument of Section 6.3 we get
\begin{eqnarray}
\int_T^{2T}F_1^4(x)dx&&=\frac{2^{5-4l}l^4B_4(\tau_l)}{(l!)^4\pi^4}\int_T^{2T}x^{
4l-1}dx\\&&+O(T^{4l+\varepsilon}y^{-1/2})
+O(T^{4l-1}V_{1,4}(y)),\nonumber
\end{eqnarray}
where
\begin{eqnarray*}
V_{1,4}(y)=\sum_{{\bf i}\in {\Bbb I}^{3} } \sum_{\stackrel{n_j\leq
y,1\leq j\leq 4}{\alpha({\bf n};{\bf i})\not= 0}}
\frac{d(n_1)d(n_2)d(n_3)d(n_4)}{(n_1n_2n_3n_4)^{3/4}}\min\left(T,\frac{T^{1/2}}{|\alpha({\bf
n};{\bf i})|}\right).\end{eqnarray*}

In \cite{Z2} the author proved that if $y\ll T^{3/4},$ then
\begin{equation}
V_{1,4}(y)\ll T^{1-3/28+\varepsilon}.
\end{equation}

From (7.10) and (7.11) we get
\begin{eqnarray}
\int_T^{2T}F_1^4(x)dx=\frac{2^{5-4l}l^4B_4(\tau_l)}{(l!)^4\pi^4}
\int_T^{2T}x^{4l-1}dx+O(T^{4l-3/28 +\varepsilon}).
\end{eqnarray}

By the argument of Section 6.4 we get
\begin{eqnarray}
\int_T^{2T}F_1^3(x)F_2(x)dx \ll T^{4l+\varepsilon}y^{-1/2}+T^{4l-1/2
}y^{1/2}{\cal L}^4 +T^{4l-1}V_{2,4}(y),
\end{eqnarray}
where
\begin{eqnarray*}
V_{2,4}(y)=\sum_{{\bf i}\in {\Bbb I}^ 3 }\sum_{y<n_1\leq 16y}
\sum_{\stackrel{n_j\leq y,2\leq j\leq 4}{\alpha({\bf n};{\bf
i})\not= 0}}
\frac{d(n_1)d(n_2)d(n_3)d(n_4)}{(n_1n_2n_3n_4)^{3/4}}\min\left(T,\frac{T^{1/2}}{|\alpha({\bf
n};{\bf i})|}\right).\end{eqnarray*}

From (7.11) we get that if $y\ll T^{3/4},$ then
\begin{equation}
V_{2,4}(y)\ll V_{1,4}(16y)\ll T^{1-3/28+\varepsilon}.
\end{equation}

From Section 7.1 we know that Lemma 6.3 holds for $y\leq
T^{77/208}.$  Taking $A=4$ in Lemma 6.3 we get
\begin{equation}
\int_T^{2T}F_2^4(x)dx\ll T^{4l+\varepsilon}y^{-77/208}.
\end{equation}
 Taking $k=4$ in (6.37) we get
\begin{equation}
\int_T^{2T}F_1^2(x)F_2^2(x)dx\ll T^{4l+\varepsilon}y^{-77/208}
\end{equation}
holds for $y\leq T^{77/208}.$

From (7.10)-(7.16) and taking $y=T^{1/3}$ we get
\begin{eqnarray*}
\int_T^{2T}F^4(x)dx=
\frac{2^{5-4l}l^4B_4(\tau_l)}{(l!)^4\pi^4}\int_T^{2T}x^{4l-1}dx
+O(T^{4l-3/28+\varepsilon}),
\end{eqnarray*}
which combining the arguments of Section 6.7  gives
\begin{eqnarray}
\int_T^{2T}R_*^4(x)dx=
\frac{2^{5-4l}l^4B_4(\tau_l)}{(l!)^4\pi^4}\int_T^{2T}x^{4l-1}dx
+O(T^{4l-3/28+\varepsilon}).
\end{eqnarray}
The formula (2.9) follows from (7.17).

\section{\bf Proofs of Theorem 5 and Theorem 6}

We shall follow Heath-Brown's argument\cite{He} to prove Theorem 5
and Theorem 6. We first quote some results  from \cite{He}. The
following Hypothesis(H), Lemma 8.1 and Lemma 8.2 are  Hypothesis(H)
 , Theorem 5 and Theorem 6 of \cite{He}, respectively.

{\bf Hypothesis(H)}: Let $M(t)$ be a real valued function, $a_1(t),
a_2(t),\cdots,$ be continuous real valued function with period 1,
and suppose there are non-zero constants $\gamma_1,\gamma_2,\cdots$
such that
$$\lim_{N\rightarrow \infty}\limsup_{T\rightarrow\infty}\frac{1}{T}\int_0^T\min\left(1,\left|M(t)-\sum_{n\leq N}a_n(\gamma_nt)\right|\right)dt=0.$$

{\bf Lemma 8.1.} Suppose $M(t)$ satisfies $(H)$ and suppose that the
constants $\gamma_i$ are linearly independently over ${\Bbb Q}.$
Suppose further
\begin{eqnarray*}
&&\int_0^1a_n(t)dt=0\ \ \ (n\in {\Bbb N}),\\
&&\sum_{n=1}^\infty\int_0^1a_n^2(t)dt<\infty,
\end{eqnarray*}
and there is a constant $\mu>1$ such that
\begin{eqnarray*}
&&\max_{t\in [0,1]}|a_n(t)|\ll n^{1-\mu},\\
&&\lim_{n\rightarrow\infty}n^\mu\int_0^1a_n^2(t)dt=\infty.
\end{eqnarray*}
Then $M(t)$ has a distribution function $f(\alpha)$ with the
properties described  in Theorem 5.

{\bf Lemma 8.2.} Suppose $M(t)$ satisfies (H) and that
\begin{eqnarray*}
\int_0^T|M(t)|^Kdt\ll T
\end{eqnarray*}
holds for some positive number $K$. Then for any real number $k\in
[0,K),$ the limit
\begin{eqnarray*}
\lim_{T\rightarrow\infty}\frac 1T\int_0^T|M(t)|^Kdt
\end{eqnarray*}
exists.

\bigskip

Suppose $T\leq x\leq 2T,$ $H=T^2, J=[({\cal L}-\log {\cal L})/2\log
2].$ Define
\begin{eqnarray*}
 M(x)&&=x^{-(2l-1/2)}R_*(x^2),\\
 a_n(x)&&=\frac{\mu^2(n)}{n^{3/4}}\sum_{r=1}^\infty
\frac{\tau_l(nr^2)}{r^{3/2}}\cos\left(2\pi rx-\frac{\pi}{4}\right),\\
 \gamma_n&&=\sqrt n.
\end{eqnarray*}
It is easy to check  that $a_n(x)$ satisfies all conditions of Lemma
8.1 for any fixed constant $3/2<\mu<7/4.$  By Proposition 6.1 we
have
\begin{eqnarray*}
&&M(x)=M_1(x)+M_2(x)+M_3(x),\\
&&M_1(x)=\frac{2^{2-l}}{(l-1)!\pi}\sum_{1\leq n\leq T{\cal
L}^{-2}}\frac{\tau_l(n)}{n^{3/4}}\cos\left(2\pi x\sqrt
n-\frac{\pi}{4}\right),\\
&&M_2(x)= \frac{2^{2-l}}{(l-1)!\pi}\sum_{T{\cal L}^{-2}< n\leq
 H^2(2^{2J+1}+1/2)}\frac{\tau_l(n;H,T)}{n^{3/4}}\cos\left(2\pi x\sqrt
 n-\frac{\pi}{4}\right),\\&&
 M_3(x)=O(T^{-1/2}G(x^2)+T^{-1/2}{\cal L}^3).
\end{eqnarray*}

It is easy to see that for any integer $N\leq T^{1/3}$ we
have(recall that $\tau_l(n;H,T)\ll d(n)$)
\begin{eqnarray*}
&&\ \ \ \ \ \ \ \ \ \ \ |M(x)-\sum_{n\leq N}a_n(\gamma_n x)|\\
&&\ll \left|\sideset{}{^{\prime}}\sum_{n\leq T{\cal
L}^{-2}}\frac{\tau_l(n)}{n^{3/4}}\cos(2\pi x\sqrt
n-\frac{\pi}{4})\right|+\sum_{n\leq N}\frac{1}{n^{3/4}}\sum_{r>\sqrt
T/\sqrt
n}\frac{d(nr^2)}{r^{3/2}} +|M_2(x)|+M_3(x)\\
&&\ll \left|\sideset{}{^{\prime}}\sum_{n\leq T{\cal
L}^{-2}}\frac{\tau_l(n)}{n^{3/4}}\cos(2\pi x\sqrt
n-\frac{\pi}{4})\right|+N^{1/2}T^{-1/4}{\cal L}^{5}+|M_2(x)|+M_3(x),
\end{eqnarray*}
where $\sideset{}{^{\prime}}\sum$ means that $n$ has a square-free
kernel great that $N.$

We have(note that $\tau_l(n)\ll d(n)$)
\begin{eqnarray*}
\int_T^{2T}\left|\sideset{}{^{\prime}}\sum_{n\leq T{\cal
L}^{-2}}\frac{\tau_l(n)}{n^{3/4}}\cos(2\pi x\sqrt
n-\frac{\pi}{4})\right|^2dx&&=T\sideset{}{^{\prime}}\sum_{n\leq
T{\cal L}^{-2}}\frac{d^2(n)}{n^{3/2}}+O(T^\varepsilon)\\
&&\ll T \sum_{n> N}\frac{d^2(n)}{n^{3/2}}+T^\varepsilon\\
&&\ll TN^{-1/2}\log^3 N+T^\varepsilon.
\end{eqnarray*}

Similar to (6.30) we have
\begin{eqnarray*}
\int_T^{2T}|M_2(x)|^2dx\ll TN^{-1/2}\log^3 N.
\end{eqnarray*}

By the trivial estimate $G(u)\ll \sqrt u(u\sim T)$ and  Lemma 6.1 we
get
\begin{eqnarray*}
&&\int_T^{2T}|M_3(x)|^2dx\ll T^{-1}\int_T^{2T}G^2(x^2)dx+{\cal L}^6\\
&&\ll \int_T^{2T}G(x^2)dx+{\cal L}^6\ll \int_{\sqrt
T}^{\sqrt{2T}}G(y)y^{-1/2}dy+{\cal L}^6 \\
&&\ll T^{1/2}H^{-1}{\cal L}+{\cal L}^6\ll  {\cal L}^6.
\end{eqnarray*}

From the above estimates and Cauchy's inequality we get
\begin{eqnarray*}
\limsup_{T\rightarrow \infty}\frac 1T\int_T^{2T}|M(x)-\sum_{n\leq
N}a_n(\gamma_n x)|^2dx\ll N^{-1/2}\log^3 N
\end{eqnarray*}
and whence Hypothesis (H) follows. From Lemma 8.1 with $\mu=5/3$ we
get Theorem 5.

Now we prove Theorem 6. According to Lemma 8.2 ,  if we can show
 that the estimate
 \begin{eqnarray}
\int_0^T x^{-k(l-1/4)}|R(x)|^{k}dx\ll T
 \end{eqnarray}
holds for  any real number $0\leq k<A_0:=262/27,$  then Theorem 6
follows directly.

Suppose $2<k<A_0$ is fixed.    Suppose $x\sim T,$ $y=T^{1/s(10)},
H=T^{5}.$ We have
$$R(x)\ll |F_1(x)|+|F_2(x)|+T^{l-1/2}G(x)+T^{l-1/2}{\cal L}^3,$$
where $F_1(x), F_2(x)$ and $G(x)$ were defined in Section 6.
 By Lemma 6.2 we get
\begin{eqnarray*} \int_T^{2T}|F_1(x)|^{10}dx\ll T^{1+10(l-1/4)},
\end{eqnarray*}
which implies that
\begin{eqnarray*}
\int_T^{2T}|F_1(x)|^kdx\ll T^{1+k(l-1/4)},
\end{eqnarray*}
By Lemma 6.3 we get
\begin{eqnarray*}
\int_T^{2T}|F_2(x)|^kdx\ll
T^{1+k(l-1/4)-(A_0-k)/2(A_0-2)s(10)+\varepsilon}\ll T^{1+k(l-1/4)}.
\end{eqnarray*}
By Lemma 6.1 we get
\begin{eqnarray*}
\int_T^{2T}T^{k(l-1/2)}G^k(x)dx\ll T^{k(l-1/2)}\
T^{(k-1)/2}\int_T^{2T} G(x)dx\ll T^{k(l-1/4)-2}.
\end{eqnarray*}
From  the above estimates we get that
\begin{eqnarray*}
\int_T^{2T}|R(x)|^kdx\ll T^{1+k(l-1/4)}
\end{eqnarray*}
and correspondingly,
\begin{eqnarray*}
\int_0^T|R(x)|^kdx\ll T^{1+k(l-1/4)},
\end{eqnarray*}
which implies (8.1) by partial summation. This completes the proof
of Theorem 6.

\noindent
Wenguang Zhai,\\
School of Mathematical Sciences,\\
Shandong Normal University,\\
Jinan, Shandong, 250014,\\
P.R.China\\
E-mail:zhaiwg@hotmail.com
\end{document}